%% file: tambyah_martin_lee_badia_2026.tex
\documentclass[oneside, reqno, 11pt, a4paper]{amsart}

\usepackage{microtype}
\AtBeginDocument{
\DeclareSymbolFont{AMSb}{U}{msb}{m}{n}
\DeclareSymbolFontAlphabet{\mathbb}{AMSb}
}

\usepackage[dvipsnames]{xcolor}
\usepackage{graphicx}
\usepackage{tikz}
\usepackage{subfigure}
\usepackage{svg}

\usepackage[a4paper, pdftex, left=2cm, top=2cm, right=2cm, bottom=2cm]{geometry}
\usepackage[final]{pdfpages}
\usepackage{changepage}

\usepackage[foot]{amsaddr}
\usepackage{caption}
\usepackage{booktabs}
\usepackage{url}

\usepackage{xr-hyper}
\usepackage{hyperref}
\hypersetup{
	breaklinks=true,
	bookmarksopen=true,
	pdftitle={GridapGeosciences},
	pdfauthor={Tamara Tambyah, Alberto Martin, David Lee, Santiago Badia},
	pdfsubject={partial differential equations},
	pdfkeywords={KEYWORDS},
	colorlinks=true,
	linkcolor=black,
	citecolor=blue,
	filecolor=black,
	urlcolor=blue
}

\usepackage[capitalise]{cleveref}
\crefname{equation}{}{}
\crefrangeformat{equation}{(#3#1#4)--(#5#2#6)}

\usepackage{csquotes}
\usepackage[backend=biber, defernumbers=false, maxbibnames=99, style=numeric-comp, isbn=false, bibencoding=utf8, safeinputenc, url=false, doi=true, giveninits=true,sorting=none,natbib=true,uniquelist=false,maxcitenames=2]{biblatex}

\usepackage{float}
\usepackage{framed}
\usepackage{verbatim}
\usepackage{fancyvrb} 
\usepackage{booktabs}
\usepackage{multirow}
\usepackage{algorithm}
\usepackage{algpseudocode}
\usepackage[inline]{enumitem}
\usepackage[printonlyused]{acronym}  
\input{acronyms.tex} 

\newtheorem{theorem}{Theorem}[section]

\newtheorem{remark}[theorem]{Remark}

\usepackage{soul}

\usepackage{amsmath, amsfonts, amssymb, amscd}
\usepackage{mathtools} 
\usepackage[only, llbracket, rrbracket]{stmaryrd}
\usepackage{upgreek}
\allowdisplaybreaks

\usepackage[mathscr]{euscript}

\usepackage{makecell}
\usepackage{listings}
\usepackage{textgreek}
\usepackage[normalem]{ulem} 

\DeclareCaptionType{listing}[Listing][List of Listings]

\usepackage{xspace}  
\input{sb-commands.tex}

\newcommand{\TS}{ T_{\surf{x}}\thesurface } 

\usepackage{lineno}

\let\oldalign\align
\let\oldendalign\endalign

\renewenvironment{align}
{\linenomathNonumbers\oldalign}
{\oldendalign\endlinenomath}

\title[GridapGeosciences.jl]{
GridapGeosciences.jl: A Julia finite element package for partial differential equations on general manifolds 
}
\date{\today}

\author[1]{Tamara A.  Tambyah$^{1}$}
\address{$^1$School of Mathematics\\Monash University\\Clayton\\Victoria 3800\\Australia }
\email{tamara.tambyah@monash.edu}

\author[2]{Alberto F. Mart\'in$^{2}$}
\email{alberto.f.martin@anu.edu.au}
\address{$^{2}$School of Computing\\Australian National University\\Canberra\\ACT 2600\\Australia }

\author[3]{David Lee$^{3}$}
\email{david.lee@bom.gov.au}
\address{ $^{3}$Bureau of Meteorology\\Melbourne\\Australia}

\author[4]{Santiago Badia$^{1\ast}$}
\email{santiago.badia@monash.edu}

\thanks{$^\ast$Corresponding author}

\AtEveryCitekey{\clearfield{eprint}} 

\begin{document}

\input{abstract.tex}

\maketitle

\input{introduction.tex}

\input{formulation.tex}

\input{implementation.tex}

\input{application.tex}

\input{conclusion.tex}

\section*{Acknowledgements}
This research is supported by the Commonwealth of Australia as represented by the Defence Science and Technology Group of the Department of Defence.
This research is also funded by the Australian Government through the Australian Research Council (project numbers DP210103092 and DP220103160). 
This research was undertaken with the assistance of resources from the Monash University and National Computational Infrastructure (NCI Australia) allocation scheme. NCI is an NCRIS-enabled capability supported by the Australian Government. 
This work is also supported by computational resources provided by the Australian Government through NCI under the National Computational Merit Allocation Scheme (NCMAS), and ANU Merit Allocation Scheme.

\subsection*{Author contribution}
All authors contributed equally to the design of the study. 
Tamara A. Tambyah drafted the article, and performed numerical simulations with the assistance of Alberto F. Mart{\'i}n.  
All authors gave approval for publication. 

\subsection*{Declaration of competing interest}
The authors declare they have no competing interests. 

\subsection*{Data availability}
This article does not contain any additional data. 
\GridapGeoscience is a free and open source library available on \href{https://github.com/gridapapps/GridapGeosciences.jl}{GitHub}.

\printbibliography

\end{document}

%% file: acronyms.tex
\acrodef{fe}[FE]{Finite Element}
\acrodefplural{fe}[FEs]{Finite Elements}
\acrodef{pde}[PDE]{Partial Differential Equation}
\acrodefplural{pde}[PDEs]{Partial Differential Equations}
\acrodef{ssprk3}[SSPRK3]{third order, strong stability preserving Runge--Kutta method}
\acrodef{api}[API]{Application Programming Interface}
\acrodef{jit}[JIT]{Just-In-Time}
\acrodef{supg}[SUPG]{streamline upwind Petrov--Galerkin}

%% file: sb-commands.tex
\DeclareMathAlphabet{\mathpzc}{OT1}{pzc}{m}{it}

\newcommand{\grad}{\nabla}
\newcommand{\Div}{\nabla \cdot}

\newcommand{\thesurface}{\mathcal{S}}
\newcommand{\surf}[1]{ \widetilde{{#1}} }
\newcommand{\surfvec}[1]{ \surf{\boldsymbol{{#1}}} }
\newcommand{\surfgrad}{ \nabla_{\thesurface} }
\newcommand{\surfdiv}{ \nabla_{\thesurface}\cdot }
\newcommand{\surflap}{ \Delta_{\thesurface} }

\newcommand{\dmap}{ \mathrm{d}\map } 
\newcommand{\V}{\mathbb{V}}

\newcommand{\en}{ {e}_{\mathrm{n}} } 
\newcommand{\re}{ \mathpzc{a} } 
\newcommand{\gravity}{ \mathpzc{g} } 
\newcommand{\coriolis}{ \mathpzc{f} } 
\newcommand{\topography}{ \mathpzc{b} } 

\newcommand{\map}{\sigma}

\newcommand{\thechart}{\mathcal{V}}
\newcommand{\chart}[1]{ {#1} }
\newcommand{\chartvec}[1]{\vec{#1} } 
\newcommand{\chartcovec}[1]{ \underline{#1} } 
\newcommand{\flux}[1]{\boldsymbol{#1}} 

\newcommand{\RT}{Raviart-Thomas\xspace}

\newcommand{\vol}{ \,\mathrm{d}\thesurface } 
\newcommand{\VOL}{ \,\mathrm{d}\chart{x} } 

\newcommand{\area}{ \,\mathrm{d}\surf{s}} 
\newcommand{\AREA}{ \,\mathrm{d}s } 

\newcommand{\trian}{\mathcal{T}} 

\newcommand{\twoD}{two-dimensional\xspace}
\newcommand{\threeD}{three-dimensional\xspace}
\newcommand{\ThreeD}{Three-dimensional\xspace}

\newcommand{\radius}{\mathscr{R}}
\newcommand{\thickness}{\mathscr{T}}

\newcommand{\Galewsky}{Galewsky\xspace}
\newcommand{\Boussinesq}{Boussinesq\xspace}

\newcommand{\GridapGeoscience}{GridapGeosciences.jl\xspace}
\newcommand{\Gridap}{Gridap.jl\xspace}
\newcommand{\gridap}{Gridap\xspace}
\newcommand{\GridapDistributed}{GridapDistributed.jl\xspace}
\newcommand{\GridapSolvers}{GridapSolvers.jl\xspace}
\newcommand{\GridapPETSc}{GridapPETSc.jl\xspace}
\newcommand{\GridapPest}{GridapP4est.jl\xspace}

\newcommand{\jl}[1]{\texttt{#1}}

\newcommand{\IntrinsicManifold}{\jl{Intrinsic\allowbreak{}Manifold()}\xspace}
\newcommand{\ExtrinsicManifold}{\jl{Extrinsic\allowbreak{}Manifold()}\xspace}
\newcommand{\ManifoldStyle}{\jl{Manifold\allowbreak{}Style}\xspace}
\newcommand{\AtlasDiscreteModel}{\jl{Atlas\allowbreak{}Discrete\allowbreak{}Model}\xspace}

\newcommand{\AtlasOctreeDistributedDiscreteModel}{\jl{Atlas\allowbreak{}Octree\allowbreak{}Distributed\allowbreak{}Discrete\allowbreak{}Model}\xspace}

\newcommand{\ExtrudedAtlasOctreeDistributedDiscreteModel}{\jl{Extruded\allowbreak{}Atlas\allowbreak{}Octree\allowbreak{}Distributed\allowbreak{}Discrete\allowbreak{}Model}\xspace}
 
\newcommand{\discretiz}{discretis\xspace}
\newcommand{\discretization}{{\discretiz}ation\xspace}

\newcommand{\discretizations}{{\discretiz}ations\xspace}

\newcommand{\RK}{Runge--Kutta\xspace}

%% file: abstract.tex
\begin{abstract}
	 
We present \GridapGeoscience, a new parallel distributed-memory Julia package for the numerical approximation of partial differential equations on general manifolds. 
Our abstract approach to generating discrete manifold domains
relies on mesh refinement in general, 
and is implemented within the \gridap finite element library,
whose high-level interface enables seamless definition of both intrinsic and extrinsic variational formulations. 
\GridapGeoscience features the cubed sphere manifold, which is used in this study to showcase the functionality of \GridapGeoscience in different  applications.
The first application involves scalar transport where we develop intrinsic formulations of common stabilisation terms, which are shown to be equivalent to well-known extrinsic stabilisation formalisations. 
Such stabilisation terms are used to develop an intrinsic formulation of the \twoD thermal shallow water equations that includes upwinded numerical fluxes. 
The final application is a linearisation of the \Boussinesq equations where we demonstrate the capacity of \GridapGeoscience to capture both surface and vertical dynamics on a \threeD manifold.

\end{abstract}

%% file: introduction.tex
\section{Introduction}

State-of-the-art \acf{fe} packages for \acfp{pde} often rely on two programming languages
\cite{Rathgeber2016,FiredrakeUserManual,Gibson2019_firedrake,Rognes2013,Logg2012_fenics,UFL2014}: 
(1) a compiled language to execute performance-critical tasks, 
and (2) an interpreted language to provide an user front-end. 
This two-language barrier is eliminated in the \gridap \ac{fe} ecosystem \cite{Verdugo2022,Badia2020,Badia2022,Manyer2024,GridapPETSc,GridapP4est},
which is developed exclusively in the Julia programming language. 
At the core of \gridap are unique abstract types \cite[Table 1]{Verdugo2022} that can be customised via multiple dispatching, 
and lazy arrays, for which the Julia 
\ac{jit} 
compiler computes entries on demand and reduces memory allocation. 
In this study, we present \GridapGeoscience, a new Julia library that extends \gridap to the numerical approximation of \acp{pde} on manifolds.
Central to \GridapGeoscience is the  expressive high-level \ac{api} of \gridap \cite{Verdugo2022,Badia2020}, which enables seamless definition of intrinsic variational forms that include general metric transformations.
Such functionality is difficult to achieve in other \ac{fe} codes that use interpreted languages to generate compiled variational forms
\cite{Rathgeber2016,FiredrakeUserManual,Gibson2019_firedrake,Rognes2013,Logg2012_fenics,UFL2014}.
High resolution simulations on the cubed sphere manifold \cite{Ronchi1996} demonstrate  
\GridapGeoscience is a high performance library capable of resolving  turbulent dynamics that occur over a variety of spatial and temporal scales.

A manifold can be represented extrinsically or intrinsically \cite{Lee2003_manifold,Frankel2004}. 
Similar to other studies \cite{Guba2014},  \GridapGeoscience
provides an extrinsic approach, where the manifold is embedded in ambient space via an exact geometrical transformation of reference shape functions.
This approach is free of geometrical errors, which arise in  studies that use a polynomial approximation of the curved geometry \cite{Rognes2013}.
In contrast, an intrinsic framework uses a smooth geometrical map to  represent the manifold via a series of
parametric spaces \cite{Lee2003_manifold,Frankel2004}.
These \textit{charts} represent local Euclidean neighbourhoods of the manifold that are connected topologically to form an \textit{atlas} \cite{Lee2003_manifold,Frankel2004}.
We develop a new approach to generating atlas triangulations that relies on a coarse representation of the atlas with a single element per chart, where mesh refinement in general yields finer charts.  
For triangles or quadrilaterals, charts are the affine mapping of a reference element,  meaning intrinsic variational forms can be numerically integrated in reference space. 
Consequently, the intrinsic approach is computationally superior to the extrinsic approach, where variational forms are computed in the ambient space of the manifold \cite{Tambyah2026_pdesphere}.
Our ability to provide both the extrinsic and intrinsic approach  
is a key novelty of  \GridapGeoscience in comparison to other  codes that exclusively provide extrinsic functionality \cite{Rognes2013}.

The \gridap ecosystem is sufficiently flexible to support atlas triangulations for manifolds in general, 
in particular the cubed sphere manifold that is featured in  \GridapGeoscience. 
The cubed sphere manifold is a spherical parametrisation often used to represent the Earth in atmospheric applications \cite{Ronchi1996,Nair2005,Staniforth2012,Mcgregor2005}. 
Such a parametrisation circumvents the `pole problem', which typically arises in other parametrisations based on latitude-longitude \cite{Staniforth2012,Putman2007}. 
Compatible \ac{fe} spaces that form a discrete de Rham complex are particularly popular in atmospheric applications 
since they preserve the necessary 2:1 ratio of velocity and pressure degrees of freedom to avoid spurious modes, and give rise to desirable conservation properties \cite{Cotter2012,Cotter2023,BauerCotter2018,Staniforth2012,Tambyah2025}. 
Such mimetic properties are satisfied for extrinsic approaches with an exact geometrical representation of the manifold, meaning many studies consider an embedding of the cubed sphere manifold in ambient space \cite{Rognes2013,Guba2014,Wimmer2020}.
This does not hold for extrinsic approaches that use $C^{0,1}$ piecewise-polynomial geometrical representations, such as FEniCS \cite{Rognes2013}, Firedrake \cite{Gibson2019_firedrake} and the LFRic dynamical core \cite{Melvin2019,Melvin2024_orography}, where the discrete velocity is not exactly tangent to the manifold \cite{Tambyah2026_pdesphere}.
The intrinsic approach in \GridapGeoscience is supported by the rigorous derivations in \citet{Tambyah2026_pdesphere}, which is our companion methods paper that uses vector calculus to systemically demonstrate mimetic and conservation properties are satisfied for an intrinsic approach.
Such findings are used in the current study to conduct intrinsic numerical experiments on the cubed sphere manifold, 
for which we explicitly define the atlas in two and three dimensions to model both surface and vertical dynamics.

This article is structured as follows. 
In \cref{sec: math formulation} we present a general intrinsic framework that is demonstrated using the well-known Poisson problem \cite[Chapter 4]{Arnold2018}. 
In \cref{sec: gridapgeociences}, we discuss the main abstractions of \gridap to manifold domains, 
and use the cubed sphere manifold to illustrate the novel design of atlas triangulations in \GridapGeoscience, which depends on mesh refinement in general.  
To demonstrate the intrinsic and extrinsic schemes are equivalent, 
we consider scalar transport on the \twoD cubed sphere manifold in \cref{sec: geo application} \cite{Lauritzen2014}.
High resolution simulations of 
the \twoD thermal shallow water equations \cite{Tambyah2025}
and the \threeD linearised \Boussinesq equations \cite{Gibson2019}
showcase the applicability of \GridapGeoscience to atmospheric case studies where turbulent non-linear flows develop, 
which extends verbatim to other geophysical applications \cite{Vallis2006}.
Potential avenues to extend the functionality in \GridapGeoscience are discussed in \cref{sec: conclusion}.

%% file: formulation.tex
\section{Mathematical formulation}

\label{sec: math formulation}

To demonstrate key principles of an intrinsic framework, we consider the following scalar Poisson problem
on an $n$-dimensional manifold  $\thesurface\subset \mathbb{R}^m$  that is closed, compact and orientable: find $ \surf{u}\in H^1(\thesurface)$ such that
\begin{align}  
	-\surflap \surf{u} = \surf{f}
	,
	\quad 
	\text{in }\thesurface,
	\label{eq: poisson}
\end{align}
where the solution $\surf{u}$ is determined up to a constant
for a given data source $\surf{f}\in L^2(\thesurface)$ that satisfies the compatibility condition $\int_{\thesurface} \surf{f} \vol = 0$, 
and $\surflap$ is the Laplace--Beltrami operator \cite{Bonito2019} defined in \cref{sec: parametric representation} below.
We now introduce the parametric representation of manifolds 
and parametric \ac{fe} methods, where \cref{eq: poisson} is used as a motivating example. 
The notation follows the companion methods paper \citet{Tambyah2026_pdesphere}, which comprehensively develops and analyses the intrinsic framework summarised below.
Readers are referred to \citet{Frankel2004} for the differential geometry background.


\subsection{Parametric representation}

\label{sec: parametric representation}

Let $n$ and $m$ be the topological and ambient dimension of the manifold $\thesurface$ where $n\leq m$. 
Typically $m=3$ in practical applications, and the topological dimension is either $n=2$ or $n=3$, which is the case hereinafter.
An \emph{atlas} $\left\{ (\thechart_k,\map_k) \right\}_{k=1}^\kappa$ of \emph{charts} $(\thechart_k,\map_k)$ parametrises the manifold,
where $\thechart_k \subset \mathbb{R}^n$ is a  Euclidean \textit{parametric domain} that is open and bounded,  $\map_k: \thechart_k \rightarrow \thesurface$ is a smooth \textit{geometrical map}, 
and $k$ is the \textit{chart identifier} \cite{Bonito2019}. 
For simplicity, 
we drop the index $k$ and consider a single chart $(\thechart,\map)$.
The following discussion applies verbatim to an atlas.

In the ambient space $\mathbb{R}^3$, let $\surf{x} = (\surf{x}^1, \surf{x}^2, \surf{x}^3)$ be a point, $\surf{f}:\thesurface\rightarrow \mathbb{R}$ be a scalar function, and $\surfvec{v} \in \TS$ be a vector in the tangent space of the manifold at $\surf{x}$, denoted $\TS$.
\ThreeD  quantities in ambient space relate to quantities in parametric space  as
\begin{align}
	\surf{x} &= \map(\chart{x}),&
	\surf{f}(\surf{x}) &= f(\chart{x}),&
	\surfvec{v}(\surf{x}) &= \dmap \, \chartvec{v}(\chart{x}),
	\label{eq: pullback}
\end{align}
where $\chart{x} = (x^1, x^2)$ are the local coordinates of $\thechart\subset \mathbb{R}^2$.
Analogous relationships hold when $\chart{x} = (x^1, x^2, x^3)$ are the local coordinates of $\thechart\subset\mathbb{R}^3$.
Hereinafter, we drop the functional dependence, and imply composition with the geometrical map.
In \cref{eq: pullback},
$\chart{f}:\thechart\rightarrow \mathbb{R}$ is the \textit{pullback} of $\surf{f}: \thesurface \rightarrow \mathbb{R}$ by the geometrical map
\cite[Chapter 9.2]{ErnGuermond_1},
and $\chartvec{v} : \thechart \rightarrow \mathbb{R}^n$ is the \emph{contravariant component field} of the tangent field $\surfvec{v}$ \cite{Tambyah2026_pdesphere}.
The affiliation between ambient and parametric quantities is implied via the same letter, where ambient quantities carry a tilde.
Further, $\dmap$ is the $m \times n$ differential of the geometrical map $\map$, also called its Jacobian,
which defines $g = \dmap^T \dmap$ as the $n \times n$ \emph{metric tensor} \cite{Frankel2004}.

The \emph{covariant component field} of $\surfvec{v}$ is $\chartcovec{v} = g \, \chartvec{v}$. 
Conversely, 
$\chartvec{v} = g^{-1} \chartcovec{v}$   \cite{Tambyah2026_pdesphere}.
The \emph{metric factor} $\sqrt{g} = (\det{g})^{1/2}$ relates integrals in the ambient and parametric spaces as
\begin{align}
	\int_{\thesurface}  \surf{f} \vol &= \int_{\thechart} \chart{f} \sqrt{g}  \VOL,
	 \label{eq: integration}
\end{align}
where $\VOL$ denotes the Lebesgue measure in the parametric space and $\vol$ is the surface measure on the manifold.

\subsubsection{Surface operators}
 \label{sec: surface operators}
For ambient scalar functions $\surf{f}: \thesurface \rightarrow  \mathbb{R}$ and 
tangent vector fields  $\surfvec{v} \in \TS$, 
the surface gradient, $\surfgrad$, and surface divergence, $\surfdiv$, and Laplace--Beltrami, $\surflap \surf{f}$,  operators are
\begin{align}
	\surfgrad \surf{f} =  \dmap \, g^{-1} \grad \chart{f},
	 \qquad 
	\surfdiv \surfvec{v} = \frac{1}{\sqrt{g}} \Div (\sqrt{g}\chartvec{v}), 
	\qquad
	\surflap \surf{f} = \surfdiv (\surfgrad \surf{f})
	,
	\label{eq: diff op}
\end{align}
where $\grad$ and $\Div$ are the Euclidean gradient and divergence operators in the parametric space, 
and $\chart{f} : \thechart \rightarrow \mathbb{R}$ and $\chartvec{v} : \thechart \rightarrow \mathbb{R}^n$ are defined via \cref{eq: pullback}
\cite{Frankel2004}. 

When $n=3$, the cross product of tangent vectors  $\surfvec{v},\surfvec{w}\in \TS$ is
\begin{align}
	\surfvec{v} \times \surfvec{w} &= \dmap\frac{\chartcovec{v} \times \chartcovec{w} }{\sqrt{g}}
	\label{eq: cross product}
\end{align}
where $\chartcovec{v}=g\chartvec{v}$, $\chartcovec{w}=g\chartvec{w}$ are covariant vectors for $\chartvec{v}$, $\chartvec{w}$  \cref{eq: pullback}. 
The surface curl and vector Laplacian operators can also be defined intrinsically when $n=3$  \cite{Frankel2004,Tambyah2026_pdesphere},
though not required for the atmospheric applications presented below in \cref{sec: geo application} of the current study.

As demonstrated in \cref{sec: tsw} below,
\twoD atmospheric applications on spherical domains  
often involve rotated vector fields, $\surfvec{v}^{\perp} $, 
the skew gradient, $\surfgrad^{\perp} \surf{f}$,
and skew divergence, $\surfgrad^{\perp} \cdot  \surfvec{v}$,  operators, defined as \cite{Tambyah2026_pdesphere}
\begin{align}
	\surfvec{v}^{\perp}  
	 =  \dmap \sqrt{g} g^{-1} R \, \chartvec{v},  
	\qquad 
	\surfgrad^{\perp} \surf{f} 
	=  \dmap \frac{\grad^\perp f}{\sqrt{g}}  ,
	\qquad 
	\surfgrad^{\perp} \cdot  \surfvec{v}
	= - \frac{1}{\sqrt{g}} \Div ( R \, \chartcovec{v} ),
	\quad 
	\text{where}
	\quad 
	R = \begin{pmatrix}
		0 &-1 \\ 1 & 0
	\end{pmatrix},
	\label{eq: skew operators}
\end{align}
is the anticlockwise rotation of \twoD vectors in the parametric space. That is, $R \chartvec{v} = (-v^2, v^1)$ for $\chartvec{v} = (v^1, v^2)$ where $\chartcovec{v}=g\chartvec{v}$,
such that
$\grad^\perp := R \, \grad = (-\partial/\partial x^2, \partial/\partial x^1)$ is the Euclidean skew gradient in $\thechart \subset \mathbb{R}^2$ \cite{Tambyah2026_pdesphere}.
Analysis of duality pairings, integration by parts, and annihilation properties for the surface operators  \cref{eq: skew operators,eq: cross product,eq: diff op} is presented in \citet{Tambyah2026_pdesphere}.
 
\subsection{Parametric finite element methods}
\label{sec: parametric FE}

To eliminate geometrical consistency errors in the discrete approximation of \acp{pde} like \cref{eq: poisson},  
we recall the parametric space is Euclidean, meaning its partition into individual elements is exact. 
Let  $\trian$  be  a  conforming triangulation of a closed subdomain  $P \subset \thechart$. 
The image of $\trian$, denoted $\surf{\trian}$, is a  conforming triangulation of the manifold $\thesurface$. 
For atlases with multiple charts $P_k \subset \thechart_k$, we collect the triangulation of each subdomain, $\mathcal{T}_k$, into the \emph{atlas triangulation} $\trian = \left\{\mathcal{T}_k \right\}_{k=1}^\kappa$, whose  image is $\surf{\trian}$, as in the single-chart case.
In general, $\int_{\thechart}$ abbreviates the sum of contributions in $\trian$, and $\V$ is a parametric \ac{fe} space  whose image is the ambient \ac{fe} space $\surf{\V}$ \cite{Tambyah2026_pdesphere}.

To construct the intrinsic \ac{fe} formulation of \cref{eq: poisson},
we consider the parametric \ac{fe} space
$\chart{\V} = \mathcal{P}_{p}(\chart{\trian})$ as the 
scalar continuous Lagrangian \ac{fe} space of piecewise polynomials of order $p$, where \cref{eq: pullback} yields the corresponding ambient \ac{fe} space $\surf{\V} = \mathcal{P}_{p}(\surf{\trian})$. 
Let $\overline{\V} = \left\{u_h \in \V : \int_{\thechart} \chart{u}_h \VOL = 0 \right\}$
where the constraint ensures a unique discrete solution.  
The intrinsic discrete system is: 
find $u_h \in \overline{\V}$ such that
\begin{align}
	\int_{\thechart} \grad \chart{v}_h \cdot (g^{-1} \grad \chart{u}_h ) \sqrt{g}  \VOL
	&= 
	\int_{\thechart} \chart{f}_h \chart{v}_h \sqrt{g} \VOL,
	\quad
	\forall \chart{v}_h \in \overline{\V}.
	\label{eq: poisson weak}
\end{align} 
To obtain \cref{eq: poisson weak}
that exactly captures the geometry of $\thesurface$ via $g^{-1}$ and $\sqrt{g}$, we use \cref{eq: integration,eq: diff op}  to expand the Laplace--Beltrami operator, and apply integration by parts  \cite{Tambyah2026_pdesphere}. 
As discussed in \citet{Tambyah2026_pdesphere},
integral-wise metric transformations can be  quantified to machine precision for a sufficient degree of numerical quadrature in the parametric space, 
meaning \cref{eq: poisson weak} is free of geometrical consistency error.
This is a key advantage of our intrinsic approach that utilises an exact triangulation of the parametric space,
in comparison to other studies that use an approximation of the curved manifold \cite{Rognes2013,FiredrakeUserManual,Mcgregor2005}.   
Other examples of intrinsic discrete systems are provided in \cref{sec: geo application} below. 

\begin{remark}
	\label{remark: fe spaces}
	The discrete systems in \cref{sec: geo application} below arise from a generic intrinsic representation of the full de Rham machinery, including the duality pairings,  integration-by-parts identities, 
	and compatible \ac{fe} spaces. 
	In particular, we consider: 
	$\V^0= \mathcal{P}_{p+1}({\trian})$ as the scalar continuous Lagrangian \ac{fe} space of piecewise polynomials of order $p+1$, 
	$\V^1 = \mathcal{RT}_p({\trian})$ as the polynomial \RT space of  order $p$ \cite{RT1977}, and 
	$\V^2 = \mathcal{P}_{p}^-({\trian})$ as the space of piecewise discontinuous polynomials of order $p$.
	The development of these global spaces on multi-chart atlases is presented in \citet{Tambyah2026_pdesphere}.
\end{remark}

\begin{remark}
	\label{remark: pullbacks}
	As in \citet{Tambyah2026_pdesphere}, vector-valued fields in $\V^1$ are represented by their \emph{flux proxy}, denoted $\flux{u} = \sqrt{g}\,\chartvec{u}$, 
	which is the parametric object represented by \RT elements. 
	The flux proxy relates to ambient tangent fields through the contravariant Piola map, $\surfvec{u} = \dmap \, \flux{u} / \sqrt{g}$, rather than \cref{eq: pullback}.
	Scalar fields of $\V^0$ and $\V^2$ are plain pullbacks  \cref{eq: pullback} \cite[Remark 4.1]{Tambyah2026_pdesphere}.
\end{remark}

%% file: implementation.tex
\section{\GridapGeoscience}
\label{sec: gridapgeociences}

We now present \GridapGeoscience, 
a new Julia library to numerically approximate \acp{pde} on  manifolds.


\subsection{Atlas triangulation}

We first present a method for generating atlas triangulations. 
Our approach  relies on
\emph{coarse atlas information} that contains two key pieces of information. 
The first is a coarse representation of the atlas with a single element per chart. This is a coarse mesh of the parametric space, denoted $\trian^0$, where the global element number is  the chart identifier $k$.
We assign chart-wise vertices to each element in $\trian^0$. 
In doing so,  the coordinates of a vertex may depend on the chart, are not necessarily unique. 
The element-wise, local-to-global node numbering is used to assign an orientation to each edge in $\trian^0$. 
Topologically, the gluing between coarse charts in $\trian^0$ can be thought of as a connected graph $G$ that is undirected, 
where each chart in the atlas is a chordless cycle in $G$ with size equal to the number of reference element vertices \cite[Chapter 1]{Chartrand2024_graphs}.
For example, a coarse mesh of quadrilaterals relates to  chordless cycles of size 4.
In \GridapGeoscience, $G$  is used implicitly to design coarse meshes, and does not have a dedicated software abstraction.

The second datum in the coarse atlas information is the geometrical map $\map_k$ of each chart, or equivalently each element in $\trian^0$.
This enables both intrinsic and extrinsic functionality, 
where a fully intrinsic framework can be obtained from coarse atlas information that includes the metric in lieu of the geometrical map.
The coarse mesh $\trian^0$ is amenable to mesh refinement in general
such that  $\trian^\ell$ is a parametric mesh with $\ell$ levels of refinement.
In $\trian^\ell$, the chart identifier $k$ of each element is equivalent to that of its root element in $\trian^0$, 
meaning the coarse atlas information is always accessible. 
This abstract methodology is applicable to manifolds in general provided  appropriate coarse atlas information and mesh refinement capability. 

\subsubsection{Cubed sphere manifold}
\label{sec: cubed sphere manifold}
Concepts of the coarse atlas information are now demonstrated using the cubed sphere manifold  \cite{Ronchi1996}. 
The atlas is composed of six charts, referred to as \textit{panels}, 
and the  geometrical map is the well-known equiangular projection \cite{Nair2005}.
Other geometrical mappings are also possible \cite{Putman2007,Wood2013,Mcgregor2005,Giraldo2003}.
We use $(\alpha, \beta)$ to denote  local \twoD chart coordinates in place of the generic local coordinates $(x^1, x^2)$ introduced in \cref{sec: parametric representation} above. 
Similarly, points in ambient space are denoted $(x, y, z)$ in lieu of $(\surf{x}^1, \surf{x}^2, \surf{x}^3)$.


For the \twoD cubed sphere manifold, 
the graph $G$ corresponds to a cube, as illustrated in \cref{fig: cubed sphere}(a).
Each  parametric space is the open square $\thechart = (-\pi/4,\pi/4)^2 $ such that the closed subdomain is $P = [-\pi/4,\pi/4]^2$ \cite{Tambyah2026_pdesphere}. 
These are identical across panels, meaning there is no global coordinate system in  parametric space.
In contrast, \threeD ambient space possess a global coordinate system, for which \cref{fig: cubed sphere}(b) shows the unique position of each panel.  
The cubic topology
illustrated in \cref{fig: cubed sphere}(c)
is free of edge permutations such that $\trian^0$ is consistently orientated\footnote{To obtain the topology in \cref{fig: cubed sphere}(c), we apply the orientation algorithm developed in \citet{Agelek2017} to the cubic topology of the CSIRO Cubic Conformal Atmospheric Model  (see Fig. 4 in \citet{Mcgregor2005}).}. 
This mesh property significantly simplifies our implementation of the \ac{fe} method, though is not a prerequisite in general.


We now derive panel-wise geometrical maps.
For Panel I that is centred on the positive $x$ axis, 
we align the local  $( \alpha, \beta)$ axes with the global  $(y,z)$ axes in ambient space (\cref{fig: cubed sphere}(b)).
This yields the relationship $y=x\tan  \alpha$ and $z = x\tan  \beta $.
Recalling
$x^2 + y^2 + z^2 = \radius^2$
for
sphere of constant radius $\radius$, 
the geometrical map for  Panel I  is 
\begin{align}
	\map_{1}( \alpha, \beta)  &=  \frac{\radius}{\rho}(1, \tan  \alpha, \tan  \beta),
	\quad 
	\text{where}
	\quad \rho^2 = 1 + \tan^2  \alpha + \tan^2  \beta,
	\label{eq: xyz panel 1}
\end{align}
with corresponding metric tensor and metric factor,
\begin{align}
	g &= \frac{\radius^2}{\rho^4 \cos^2 \alpha \cos^2 \beta}
	\begin{pmatrix}
		1+\tan^2 \alpha & -\tan  \alpha \tan  \beta \\
		-\tan \alpha\tan \beta  & 1+\tan^2 \beta
	\end{pmatrix},
	&
	\sqrt{g} &= \frac{\radius^2}{\rho^3 \cos^2  \alpha \cos^2  \beta } .
	\label{eq: riemanniam metric}
\end{align}
Accounting for the position and local orientation of each panel with respect to global ambient space yields the panel-wise geometrical mapping  as
\begin{align}
	\text{Panel I:} ~
	\map_{1}( \alpha, \beta) &=  \frac{\radius}{\rho}(1,\tan  \alpha,\tan  \beta) ,
	&
	\text{Panel IV:} ~
	\map_4( \alpha, \beta) &=  \frac{\radius}{\rho}(-1,\tan  \beta,\tan  \alpha) , \nonumber
	\\
	\text{Panel II:} ~
	\map_2( \alpha, \beta) &=  \frac{\radius}{\rho}(-\tan \beta, \tan  \alpha,1) ,
	&
	\text{Panel V:} ~
	\map_5( \alpha, \beta) &=  \frac{\radius}{\rho}(-\tan \alpha,\tan  \beta,-1) ,
	\label{eq: panel wise projection}
	\\
	\text{Panel III:} ~
	\map_3( \alpha, \beta) &=  \frac{\radius}{\rho}(-\tan  \alpha,1,\tan  \beta) ,
	&
	\text{Panel VI:} ~
	\map_6( \alpha, \beta) &=  \frac{\radius}{\rho}(-\tan \beta,-1,\tan  \alpha). \nonumber
\end{align} 
The Riemannian metric and measure of each geometrical map in \cref{eq: panel wise projection} is  equivalent to \cref{eq: riemanniam metric}
since $\mathrm{d}\map_k = R_k \, \mathrm{d}\map_1$,
where $R_k\in \mathbb{R}^{3\times3} $ is a rotation matrix for $k=1,2\dots,6$.

\input{cubed_sphere.tex}

\subsubsection{Extension to topological dimension 3}

\label{sec: 3D CS}

A \threeD cubed sphere manifold is also available in \GridapGeoscience. 
Such a manifold has variable radius, topological dimension $n=3$, and local  coordinates $( \alpha, \beta, \gamma)$.
The extrusion in 
\cref{fig: cubed sphere}(c) illustrates $\trian^0$, 
where the orientation of  $( \alpha, \beta)$ in three dimensions corresponds to that in two dimensions, and $ \gamma$ points in the outward radial direction.

There are two main ways to extend the \twoD geometrical maps  \cref{eq: panel wise projection} to $n=3$. 
The first approach is to modify \cref{eq: panel wise projection} to include $\radius( \gamma)$.
The second method is to extrude \cref{eq: panel wise projection} in the radial direction as
\begin{align}
	\map_k( \alpha, \beta, \gamma)
	&=  \map_k( \alpha, \beta) 
	+ \topography(\alpha,\beta)
	+ 
	 \gamma [\thickness-\topography( \alpha, \beta)]  \frac{\map_k( \alpha, \beta)}{|| \map_k( \alpha, \beta) ||}  ,
	&
	k&=1,2,\dots,6,
	\label{eq: 3D map}
\end{align}
where  $( \alpha, \beta) \in [-\pi/4,\pi/4]^2$ and 
$ \gamma \in [0,1]$ such that $\gamma=0,1$ relates to the radial boundary of  $\thesurface$, and $\thickness \ll \radius$ is the thickness. 
In \cref{eq: 3D map},  $\topography(\alpha,\beta)$ is a topography profile that displaces the bottom shell at $\gamma=0$, while the top shell at $\gamma=1$ is flat. 
Setting $\topography = 0$ in \cref{eq: 3D map} yields a  \threeD atmospheric shell with radial coordinate $\radius + \thickness \gamma$ such that \cref{eq: panel wise projection} is recovered for $\gamma = 0$. 
To represent a domain with variable orography \cite{Tambyah2026_pdesphere}, we use $\topography(\alpha,\beta)$ such that the radial coordinate is $ \radius + \topography(\alpha,\beta) + \gamma(\thickness - \topography(\alpha,\beta))$.


\subsection{Implementation}

The user-level \AtlasDiscreteModel  data structure in \GridapGeoscience
represents an atlas triangulation with $\ell$ levels of refinement that originates from  coarse atlas information.  
This data structure possesses a manifold style trait, of type \ManifoldStyle, which can have one of two possible
values: either  \IntrinsicManifold or  \ExtrinsicManifold.
This trait is used under-the-hood to determine the element-wise mapping  from the \textit{element reference domain} to the \textit{physical domain} where the weak formulation is posed \cite{Badia2018fempar}.
For  \IntrinsicManifold, the physical domain represents the parametric space of the manifold, where elements $\chart{K}=\phi_K(\widehat{K})$ are the image of an affine mapping, $\phi_K$, of reference elements,  $\widehat{K}$. 
In contrast, 
\ExtrinsicManifold means the physical domain relates to the ambient space of the manifold, where  $\surf{K} = \map_k \circ \phi_K(\widehat{K})$ is an ambient element that originates from chart $k$.
Customising the physical domain of  the \AtlasDiscreteModel in this way enables both intrinsic and extrinsic functionality in \GridapGeoscience.

Users of \GridapGeoscience define variational formulations in the physical domain via \gridap's high-level 
\ac{api}. 
The automatic engine of \gridap handles the transformation of 
integrals from the physical to the reference domain under-the-hood \cite{Badia2020,Verdugo2022}. 
Importantly, this automatic engine is flexible enough to support tailored element-wise mappings that are not necessarily polynomial, but defined as the composition of polynomial maps and analytical closed-form maps. 
Ultimately element-wise mappings are extracted from the \AtlasDiscreteModel.
For \IntrinsicManifold, we assume the users of \GridapGeoscience to be informed enough to include metric-related transformations in their intrinsic variational forms.
For \ExtrinsicManifold, variational formulations are directly written in ambient space without explicit metric-related transformations.
The extensible core of \gridap is used to implement the chart-wise geometrical map $\map_k$ via
the \jl{Field} interface such that $\dmap_k$ is the associated \jl{FieldGradient}.
These abstractions showcase our ability to customise \gridap to atlas triangulations without modifying the library's core functionality. 
In large part, these extensions are supported since \gridap 
provides a modular low-level interface based on multiple dispatch  and an expressive high-level \ac{api}, 
all exclusively in the Julia programming language. 
Similar functionality may be challenging to achieve in other \ac{fe} frameworks that incur a two-language barrier 
\cite{Rathgeber2016,FiredrakeUserManual,Gibson2019_firedrake,Rognes2013,Logg2012_fenics,UFL2014}.

\subsubsection{Example usage}
\label{sec: example usage}

\cref{lst: laplace beltrami} shows a \GridapGeoscience driver that solves the Poisson problem \cref{eq: poisson} using the intrinsic \ac{fe} formulation \cref{eq: poisson weak} with a manufactured solution.
The main function (lines \ref{serial-lstMainfuncStart}--\ref{serial-lstMainfuncEnd})
is executed in line \ref{serial-lstExecCall} for an \AtlasDiscreteModel  with $\ell=3$ levels of refinement.
Such a discrete atlas is initialised in lines \ref{serial-lstAtlasStart}--\ref{serial-lstAtlasEnd},
where \jl{CubedSphereMesh(1.0)} represents the coarse mesh $\trian^0$ for the \twoD cubed sphere manifold with $\radius=1$ and geometrical map \cref{eq: panel wise projection}.
Modifying  to \jl{CubedSphereWithThicknessMesh($\radius$,$\thickness$)} yields the \threeD  cubed sphere manifold with the extruded geometrical map \cref{eq: 3D map}.
The style trait \IntrinsicManifold (line \ref{serial-lstAtlasEnd})
reflects the intrinsic nature of the weak formulation \cref{eq: poisson weak}.
\begin{listing}[h!]
	\centering 
	\includegraphics[width=0.96\textwidth]{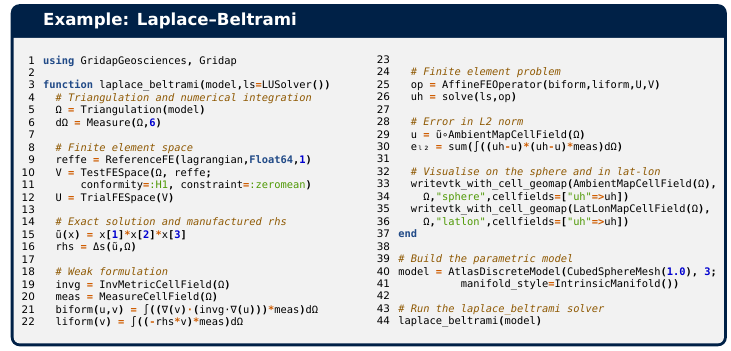}
	\caption{Driver for the intrinsic formulation \cref{eq: poisson weak} for the	\twoD cubed sphere manifold with $\radius=1$ and $\ell=3$.
	}
	\label{lst: laplace beltrami}
\end{listing}

In the main function, the triangulation (line \ref{serial-lstTriangulation}),
numerical integration (line \ref{serial-lstNumIntegration}),
\ac{fe} spaces (lines \ref{serial-lstFespaceStart}--\ref{serial-lstFespaceEnd}),
and \ac{fe} problem (lines \ref{serial-lstFeproblemStart}--\ref{serial-lstFeproblemEnd}) follow
the \gridap interface; refer to the \gridap tutorials \cite{GridapTutorials}.
To manufacture $\surf{f} = - \surflap{\surf{u}_{\mathrm{ex}}}$,
we define the Julia function $\surf{\jl{u}}$ in line \ref{serial-lstUExact} that corresponds to the ambient analytic function $\surf{u}_{\mathrm{ex}}=xyz$.
The manufactured \jl{rhs} (line \ref{serial-lstRhs}) is used in conjunction with
\jl{CellField} representations of $g^{-1}$ (line \ref{serial-lstInvg}) and $\sqrt{g}$ (line \ref{serial-lstMeas})
to define the intrinsic formulation in lines \ref{serial-lstWeakformStart}--\ref{serial-lstWeakformEnd},
where the direct correspondence to \cref{eq: poisson weak} is a key advantage of the expressive \ac{api} of \gridap.
The exact intrinsic solution \jl{u} in line \ref{serial-lstComposeU} is the pullback of $\surf{\jl{u}}$ via \cref{eq: pullback},
which is further used to compute the error of the intrinsic numerical solution \jl{u}$_\jl{h}$ in the $L^2$ norm  (line \ref{serial-lstL2error}).
Line \ref{serial-lstRhs} is the Laplace--Beltrami operator \cref{eq: diff op}. Other surface operators introduced in \cref{sec: surface operators} above are also available.

The \jl{CellField} maps   returned by the \jl{AmbientMapCellField} and the \jl{LatLonMapCellField} calls
represent the geometrical map of the cubed sphere \cref{eq: panel wise projection}, and the standard conversion from \threeD coordinates to latitude longitude coordinates, respectively. 
Such \jl{CellField} maps are only used in \cref{lst: laplace beltrami} to convert ambient input data to parametric input data (line \ref{serial-lstComposeU}), and to generate data files in VTK format to visualise,  via external software, intrinsic data in the ambient space  (lines \ref{serial-lstVisualiseStart}--\ref{serial-lstVisualiseEnd}).  	
No \ac{fe} operation is performed with these \jl{CellField} maps.

\subsection{Parallelisation}
\GridapGeoscience also provides a fully parallel distributed-memory implementation of intrinsic and extrinsic \ac{fe} methods on manifolds. 
To this end, \GridapGeoscience, on one hand,  relies  on \GridapDistributed \cite{Badia2022}, 
the parallel distributed-memory extension of \Gridap \cite{Badia2020,Verdugo2022}. 
On the other hand, scalable atlas triangulations are implemented using \GridapPest \cite{GridapP4est}.
This latter package provides forests-of-quadtrees ($n=2$) and 
forests-of-octrees ($n=3$) endowed with the so-called Morton (also known as $Z$-shaped) 
space filling curve 
\cite{Burstedde2011,Isaac2015,Isaac2015_thesis} for efficient storage, and scalable mesh refinement and partitioning.
With this approach,  the domain is represented using two levels: 
(1) a macro level that contains the coarse mesh, 
and (2) a micro level where each element is the root of an adaptive tree that can be recursively refined. 
For example, the macro coarse mesh in two dimensions is a conforming quadrilateral mesh where isotropic 1:4 refinement yields an adaptive tree-based mesh  at the micro level \cite{Burstedde2011,Isaac2015}. 
Such \textit{quadtrees} constitute the forest-of-quadtrees data structure.
Similarly in three dimensions, isotropic 1:8 refinement is applied to a hexahedral coarse mesh to obtain a forest-of-octrees \cite{Burstedde2011,Isaac2015}. 
\GridapPest also supports refinement of anisotropic \threeD domains, where the horizontal scale is much larger than the vertical scale. For such domains, a forest-of-quadtrees is extruded into layers that  represent a vertical column of elements in three dimensions. 
Consequently, the micro level mesh is obtained by refining layers of quadrilaterals in parallel \cite[Chapter 2.3]{Isaac2015_thesis}. Ultimately, the algorithms on which \GridapPest relies are based on the p4est library \cite{Burstedde2011,Isaac2015_thesis} and its p6est extension \cite{Muller2019_p6est,Isaac2015_thesis}.

In \GridapGeoscience, the \AtlasOctreeDistributedDiscreteModel \ data structure is a fully parallel, 
tree-based atlas triangulation with $\ell$ levels of isotropic refinement
and underlying coarse atlas information. 
When $n=2$, the macro coarse mesh $\trian^0$ is a forest-of-quadtrees data structure 
such that $\trian^\ell$ is the successive refinement of quadrilaterals at the micro level.
The associated forest-of-octrees data structure is applicable when $n=3$.
In geophysical applications,  $n=3$  typically represents an anisotropic domain.
Thus, in \GridapGeoscience, we also provide the \ExtrudedAtlasOctreeDistributedDiscreteModel \ data structure
to enable anisotropic refinement of a coarse quadrilateral mesh. 
Similar to the serial atlases, the physical domain of the distributed atlases is determined by the style trait, 
\IntrinsicManifold or \ExtrinsicManifold, which enables  both intrinsic and extrinsic functionality.


\subsubsection{Example usage in parallel}

\begin{figure}[h!]
	\centering
	\includegraphics{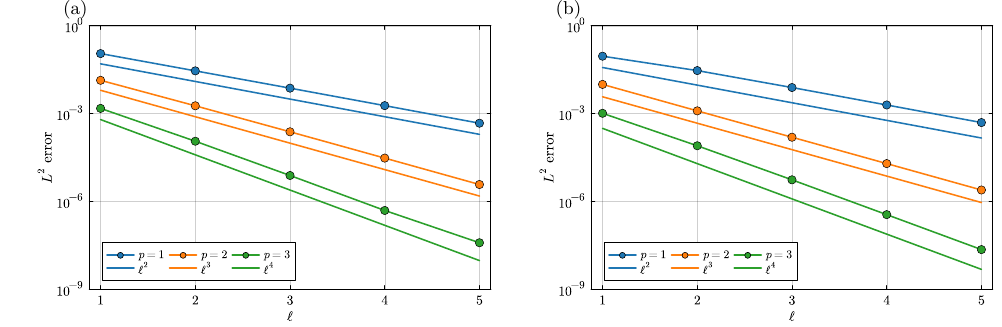}
	\caption{
		Convergence of the Poisson equation on the cubed sphere manifold showing the $L^2$ norm between the exact and numerical solutions for (a) $n=2$ and (b) $n=3$, with $6\times 2^{n\ell}$ spatial elements. 
		Parameters: $\surf{u}_{\mathrm{ex}}=xyz$, $\radius=1$, $\thickness = 0.19$.
	}
	\label{fig: laplace}
\end{figure}

\cref{lst: laplace beltrami distributed} shows a  distributed driver 
that executes the main function of \cref{lst: laplace beltrami} in multiple parallel processes.
The \AtlasOctreeDistributedDiscreteModel built in lines \ref{dist-lstAtlasStart}--\ref{dist-lstAtlasEnd} is a fully parallel atlas triangulation that is initialised with the distributed rank indices of the parallel processes, coarse atlas information associated to the \twoD cubed sphere manifold, and
$\ell=3$ levels of refinement.
Modifying   line \ref{dist-lstAtlasStart} to  \ExtrudedAtlasOctreeDistributedDiscreteModel
and line \ref{dist-lstAtlasMesh} to \jl{ExtrudedCubedSphereWithThicknessMesh($\radius$,$\thickness$)}
yields the \threeD  cubed sphere manifold, described in \cref{sec: 3D CS} above, that is amenable to anisotropic refinement.
We rely on  \GridapSolvers \cite{Manyer2024} and \GridapPETSc \cite{GridapPETSc} in general to reduce the computational bottleneck associated to solving sparse linear systems related to \ac{fe} \discretizations.
The solver in lines \ref{dist-lstSolverStart}--\ref{dist-lstSolverEnd} is provided via \GridapSolvers;
see \gridap tutorials for usage details of  \GridapPETSc \cite{GridapTutorials}. 
\Cref{fig: laplace} shows that applying \cref{lst: laplace beltrami distributed}  to a sequence of $h$-refined meshes yields the expected $hp$-convergence rates on two- and \threeD cubed sphere manifolds.
For the \threeD case, the right-hand side of \cref{eq: poisson weak} includes the flux term that arises from integration by parts on the radial boundary, 
where boundary data for the manufactured solution is prepared following Remark 4.2 in \citet{Tambyah2026_pdesphere}. 
\begin{listing}[h!]
	\centering 
	\includegraphics[width=0.96\textwidth]{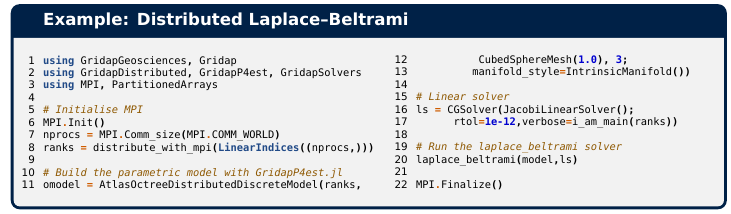}
	\caption{Parallel distributed-memory  version of \cref{lst: laplace beltrami}, 
	to be executed in multiple parallel processes via \jl{mpiexec} \cite{Byrne2021}. }
	\label{lst: laplace beltrami distributed}
\end{listing}

%% file: cubed_sphere.tex
\begin{figure}[h!]
	\centering 
	\includegraphics[width=0.9\textwidth]{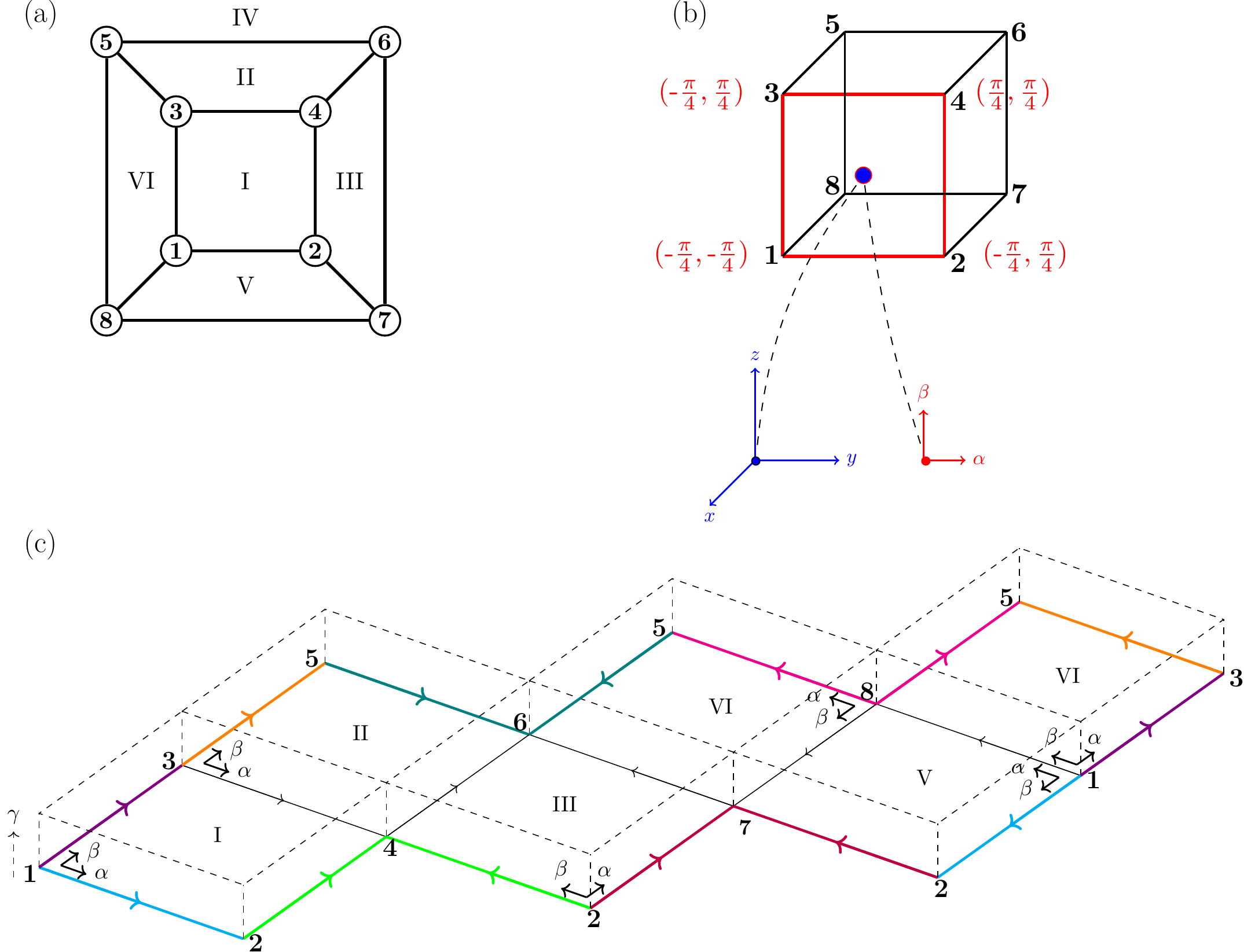}
	\caption{The cubed sphere manifold. 
		(a) is the graph $G$  that 
		relates to the  cube in (b), 
		where the local coordinates of Panel I, shown in red, align with the global $(y,z)$ coordinates shown in blue. 
		(c) illustrates the orientation of the edges in $G$ with arrows in two dimensions, 
		where colours indicate edge periodicity and the dashed lines indicate the extrusion to three dimensions. 
		The node numbering is consistent between the three sub-figures, and I--VI denote panel number.  
	}
	\label{fig: cubed sphere}
\end{figure}

%% file: application.tex
\section{Application}
\label{sec: geo application}
 
The features of \GridapGeoscience are now demonstrated on two- and \threeD cubed sphere manifolds for different atmospheric applications with \ac{fe} spaces and pullbacks as in \cref{remark: fe spaces,remark: pullbacks} above, respectively.
In all simulations, 
the number of spatial elements is $N = 6 \times 2^{n\ell}$ such that  $\ell \rightarrow L$ relates to $h\rightarrow 0$, 
$\lambda = \arctan(y/x)$ is longitude,
and $\theta = \arcsin(z/r)$ is latitude, where $r^2 = x^2 + y^2 + z^2$.
Tests are conducted on the Gadi supercomputer at NCI Australia, 
and drivers are available in the \GridapGeoscience tutorials \cite{GridapGeoscienceTutorials}.


\subsection{Scalar transport}
\label{sec: scalar transport}

The first application is transport of a scalar $\surf{u}: \thesurface \rightarrow \mathbb{R}$ by a known tangent velocity field  $\surfvec{w} \in \TS $ in the absence of sources or sinks, described by the advection equation 
\begin{align}
	\partial_t \surf{u}  + \surfdiv (\surfvec{w} \surf{u} ) = 0 , 
	\quad 
	\text{in }\thesurface \times (0,t],
	\label{eq: advection flux}
\end{align}
where $\thesurface$ is a \twoD manifold. 
As \ac{fe} formulations of such transport processes often incur spurious oscillations at the semi-discrete level,  
we consider two stabilisation schemes: 
(1) \ac{supg} method 
for continuous polynomial solutions, and (2) upwinding stabilisation for discontinuous polynomial solutions. 
We choose to formulate the \ac{supg} method for $\surfdiv \surfvec{w} = 0$ and the upwinding scheme for $\surfdiv \surfvec{w} \neq  0$. 
In both cases, we develop the intrinsic formulation, 
and compare to the well-known extrinsic formulation.

\subsubsection{Streamline upwind Petrov--Galerkin method}
For divergence free velocity fields, the \ac{fe} formulation of \cref{eq: advection flux} in the ambient space is:
find $\surf{u}_h \in \surf{\V}$ such that
\begin{align}
	\label{eq: advection supg extrinsic}
	\surf{a}(\surf{u}_h,\surf{v}_h) &+ \tau \surf{a}(\surf{u}_h,\surfvec{w}\cdot \surfgrad\surf{v}_h)  = 0 , 
	\quad 
	\text{where}, 
	\quad 
	\surf{a}(\surf{u}_h,\surf{v}_h) 
	=\int_{\thesurface} \partial_t \surf{u}_h
	\surf{v}_h \vol  
	+ \int_{\thesurface} (\surfvec{w} \cdot \surfgrad \surf{u}_h )\surf{v}_h
\vol , 
\end{align}
for all $\surf{v}_h \in \surf{\V} $, where
$\tau \in \mathbb{R}$ is the stabilisation parameter. 
The stabilisation term results from weighting test functions toward the upstream direction \cite{BrooksHughes1982}. 
Using \cref{eq: diff op,eq: pullback}, 
and recalling that $\partial_t \surf{u}_h  = \partial_t \chart{u}_h$
since the geometrical map  is a space-only map,
yields the corresponding intrinsic formulation:
find $\chart{u}_h \in \V$ such that 
\begin{align}	
	\label{eq: advection supg}
	a(\chart{u}_h,\chart{v}_h) &+ \tau a(\chart{u}_h,\chartvec{w}\cdot \grad \chart{v}_h)  = 0 ,  
	~
	\text{where}, 
	~
	a(\chart{u}_h,\chart{v}_h) 
	= 	\int_{\thechart} \partial_t \chart{u}_h
	\chart{v}_h \sqrt{g}\VOL 
	+ \int_{\thechart} (\chartvec{w} \cdot \grad \chart{u}_h )\chart{v}_h
	\sqrt{g} \VOL , 
\end{align}
for all $ \chart{v}_h \in \V$. 
The integrals in \cref{eq: advection supg} are instances of the metric-weighted pairings derived systematically in \citet{Tambyah2026_pdesphere}.
\cref{lst: supg} is the implementation of \cref{eq: advection supg}
for a given initial scalar field and time-independent velocity field.
The driver follows  \cref{lst: laplace beltrami} for the most part, with the following key differences. 
In line \ref{supg-lstVelocity}, the intrinsic velocity field is obtained via the contravariant mapping \cref{eq: pullback}.
The semi-discrete system (lines \ref{supg-lstSemidiscreteStart}--\ref{supg-lstSemidiscreteEnd}) is integrated in time using a \RK method (lines \ref{supg-lstRkStart}--\ref{supg-lstRkEnd});
refer to \gridap Tutorial 17 for details of the transient interface \cite{GridapTutorials}. 
\cref{lst: supg} supports a flexible choice of solver and Butcher tableau, and extends to parallel in a similar fashion to \cref{lst: laplace beltrami distributed}. 
\begin{listing}[h!]
	\centering 
	\includegraphics[width=0.96\textwidth]{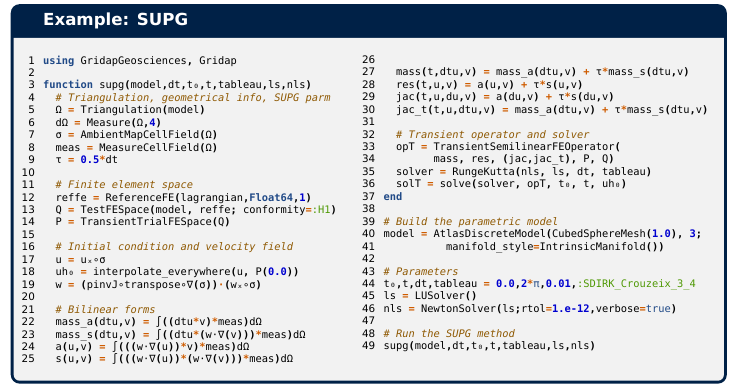}
	\caption{\GridapGeoscience driver for the intrinsic \ac{supg} method  \cref{eq: advection supg} where \jl{u$_{\jl{x}}$} and ${\textbf{\jl{w}}}_{\jl{x}}$
	are user defined Julia functions in the ambient space.
	}
	\label{lst: supg}
\end{listing}

\subsubsection{Upwinding stabilisation}
\label{sec: upwinding}
For divergent vector fields, \cref{eq: advection flux} is a conservation law, which is commonly discretised using the discontinuous Galerkin method with upwinding stabilisation. 
Let $\mathcal{E}_0$ be the set of interior edges in the atlas triangulation whose image is $\surf{\mathcal{E}_0}$. 
On each edge $e \in \mathcal{E}_0$,
the Euclidean outward unit normal, $\chartcovec{n}$, is regarded as a covariant field,
and integration over $e$ is associated to the Euclidean measure $\AREA$.
For each ambient edge $\surf{e} \in \surf{\mathcal{E}}_0 $,
the  outward pointing \emph{conormal} $\surfvec{n}$ is normal to $\surf{e}$, defined as
\begin{align}
	\surfvec{n} = \frac{ \dmap \, g^{-1} \chartcovec{n} }{ || \dmap \, g^{-1} \chartcovec{n} || } .
	\label{eq: push normal}
\end{align}
This conormal is a unit vector that lies in the tangent space of the manifold. 
Integration over $\surf{e}$ is associated to the measure $\area = \sqrt{g}||\dmap \, g^{-1} \chartcovec{n}||\AREA$  \cite[Chapter 4]{Bonet2008}.
That is,
\begin{align}
	\int_{\surf{\mathcal{E}}_0}  \surf{f} \area
	&= \int_{\chart{\mathcal{E}_0}} \chart{f}
	\sqrt{g} ||\dmap \, g^{-1} \chartcovec{n} || \AREA ,
	\label{eq: integration edge}
\end{align}
where $\surf{f}: \thesurface \rightarrow \mathbb{R}$ and $\chart{f}: \thechart \rightarrow \mathbb{R}$ are scalar functions related via \cref{eq: pullback}.

In the ambient space, there exists a global coordinate system. 
Thus, $\surfvec{n}^+ = - \surfvec{n}^-$ for  neighbouring spatial elements $\surf{K}^\pm \in \surf{\trian}$ that share a common edge $\surf{e}$, and standard jump and average operators that involve normal vectors apply \cite{Brezzi2004}. 
These conventional jump and average operators are also applicable to flat manifolds and curved manifolds parameterised by a single chart. 
For an atlas comprised of multiple charts that have different coordinate systems, such as the cubed sphere manifold in  \cref{sec: cubed sphere manifold}, a global parametric coordinate system cannot be defined. 
Consequently, $\chartcovec{n}^+ \neq - \chartcovec{n}^-$ for neighbouring spatial elements in parametric space $\chart{K}^\pm \in \chart{\trian}$ that share a common edge $\chart{e}$ on the interface of charts.
Thus, the following weak formulations are constructed using $\langle \surf{a}\rangle  = \langle a \rangle = a^+ - a^-$ where $a$ is any scalar quantity,
as opposed to jump and average operators that involve normal vectors.

In the ambient space, the semi-discrete approximation of \cref{eq: advection flux} that uses an upwinded jump penalty term  is \cite{Brezzi2004}: 
find $\surf{u}_h \in \surf{\V^2}$ such that 
\begin{subequations}
	\label{eq: weak advection extrinisc}
	\begin{align}
		\surf{a}(\chart{u}_h,\surf{v}_h) & + \surf{s}(\surf{u}_h,\surf{v}_h) = 0  ,
		\quad \text{where},
		\\
		\surf{a}(\surf{u}_h,\surf{v}_h) 
		&= 
		\int_{\thesurface} \partial_t \surf{u}_h \surf{v}_h \vol
		- \int_{\thesurface} \surf{u}_h\surfvec{w} \cdot \surfgrad \surf{v}_h \vol
		+
		 \int_{\surf{\mathcal{E}}_0} 
		\frac{1}{2} 
		\langle \surf{u}_h \surfvec{w} \cdot \surfvec{n}
		\rangle
		\langle \surf{v}_h 
		\rangle
		\area , 
		\\
		\surf{s}(\surf{u}_h,\surf{v}_h) 
		&=   
		  \int_{\surf{\mathcal{E}}_0} \frac{1}{2}  
		  |\surfvec{w}\cdot \surfvec{n}^+|
		  \langle  \surf{u}_h \rangle
		    \langle \surf{v}_h  \rangle 
		  \area
		  ,  
	\end{align} 
\end{subequations}%
for all $\surf{v}_h \in \surf{\V^2}$. 
Following the systematic derivations in \citet{Tambyah2026_pdesphere}, we use \cref{eq: push normal,eq: diff op,eq: pullback,eq: integration,eq: integration edge} to relate ambient and parametric quantities, which yields the intrinsic formulation:
find $\chart{u}_h \in \V^2$ such that 
\begin{subequations}
		\label{eq: weak advection chart}
	\begin{align}
		a(\chart{u}_h,\chart{v}_h)   &+ s(\chart{u}_h,\chart{v}_h) = 0,  \quad \text{where},
		\\
		a(\chart{u}_h,\chart{v}_h) 
		&= 
		\int_{\thechart} \partial_t \chart{u}_h \chart{v}_h \sqrt{g} \VOL
		- \int_{\thechart} \chart{u}_h\chartvec{w} \cdot \grad \chart{v}_h  \sqrt{g}\VOL 
		+ \int_{\chart{\mathcal{E}}_0} 
		\frac{1}{2}\langle  \chart{u}_h \chartvec{w} \cdot \chartcovec{n}
		\rangle 
		\langle \chart{v}_h \rangle \sqrt{g}\AREA,
		\\
		s(\chart{u}_h,\chart{v}_h) &= \int_{\chart{\mathcal{E}}_0} \frac{1}{2}  
		|\chartvec{w}\cdot \chartcovec{n}^+|
		\langle \chart{u}_h \rangle 
		\langle \chart{v}_h \rangle \sqrt{g}  \AREA , 
		\label{eq: weak advection chart stablisation}
	\end{align}
\end{subequations}%
for all $\chart{v}_h \in \V^2$. 
\cref{lst: upwinding} is the driver for 
\cref{eq: weak advection chart}, where lines \ref{upw-lstBraketStart}--\ref{upw-lstBraketEnd} are the function that represents $\langle a \rangle$,
which is used in semi-discrete system in lines \ref{upw-lstSemidiscreteStart}--\ref{upw-lstSemidiscreteEnd}.
Otherwise, the main function in \cref{lst: upwinding} closely follows that of \cref{lst: supg}.
\begin{listing}[h!]
	\centering 
	\includegraphics[width=0.96\textwidth]{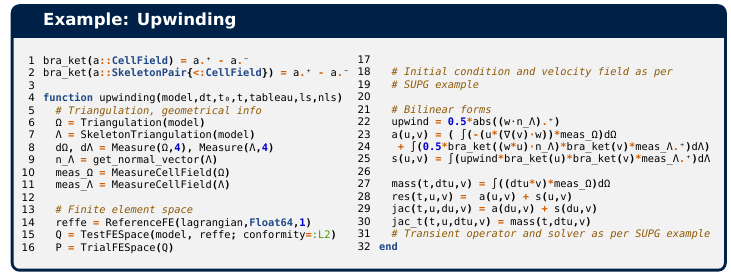}
	\caption{
	\GridapGeoscience driver for the intrinsic upwinding scheme \cref{eq: weak advection chart}, 
	where the \jl{bra\_ket} function represents $\langle a \rangle$.
	This driver is executed as per  \cref{lst: supg}.
	}
	\label{lst: upwinding}
\end{listing}

\subsubsection{Results}
\begin{figure}[h!]
	\centering
	\includegraphics{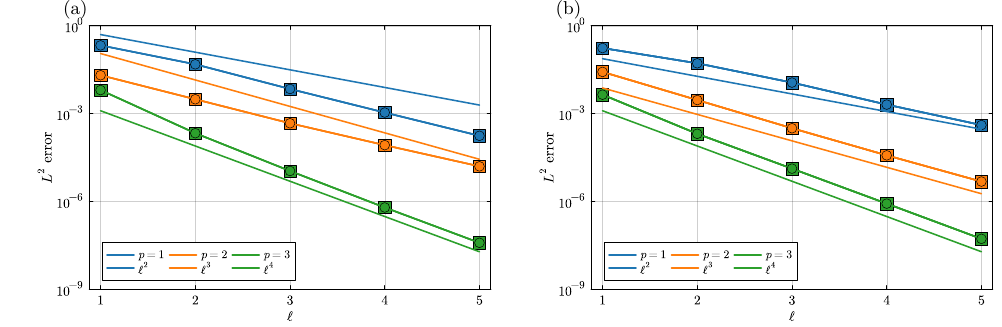}
	\caption{
	Convergence of the solid body rotation test case showing the $L^2$ norm between the initial and final solutions
	for (a) the \ac{supg} method and (b) the upwinded scheme, 
	with $6\times 2^{2\ell}$ spatial elements. 
	Circles and squares indicate the parametric and ambient solution, respectively. 
	Parameters: $\radius = 1$, $\mathrm{CFL} = 0.1$.
	}
	\label{fig: advection}
\end{figure}

A solid body rotation test case is used to  compare the intrinsic formulations \cref{eq: advection supg}, \cref{eq: weak advection chart}
and the extrinsic formulations \cref{eq: advection supg extrinsic}, \cref{eq: weak advection extrinisc}. 
The initial scalar field and time-independent velocity field is
\begin{align}
	\surf{u}(0) &= \exp(-( y^2 + z^2 ) ) ,
	&
	\surfvec{w} &= (-y,x,0) .
	\label{eq: advection IC}
\end{align}
Convergence is assessed using the $L^2$ norm between the initial and final solution at $t=2\pi$. 
To minimise errors due to temporal integration, 
the semi-discrete systems 
\cref{eq: advection supg,eq: weak advection chart,eq: advection supg extrinsic,eq: weak advection extrinisc}
are integrated using  Crouzeix's three stage, fourth order, diagonally implicit \RK method with time step $\Delta t = h \mathrm{CFL}$. 
At each \RK stage, the \ac{fe} linear system is solved via the UMFPACK sparse direct solver in Julia. 
\Cref{fig: advection} shows the expected $hp$-convergence rates for both the \ac{supg} method with $\tau = \Delta t/2$, and the upwinded scheme. 
For each line in \cref{fig: advection}, $h$ and $\Delta t$ reduce simultaneously, and the extrinsic and intrinsic solutions are equivalent. 
Hereinafter, we exclusively use an intrinsic approach. 

\begin{figure}[h!]
	\centering
	\includegraphics{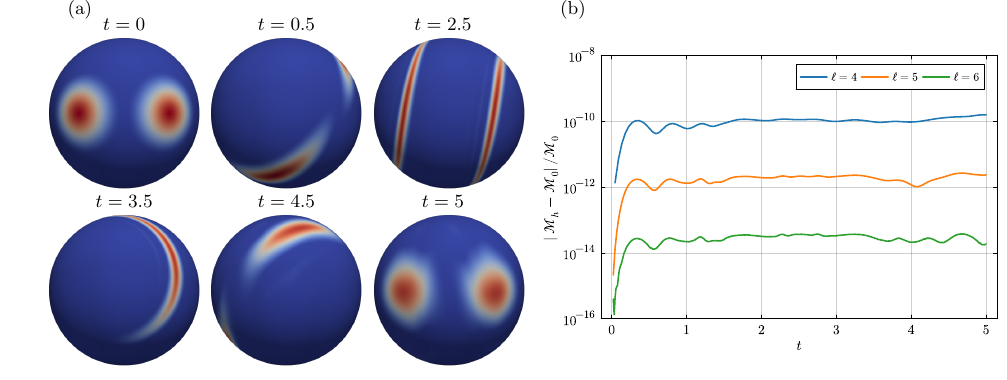}
	\caption{Advection of two cosine bells where (a) shows the scalar field at various times for $\ell=6$, 
	and (b) shows mass conservation error for various refinement levels, where time is unit-less. 
	Parameters: $\radius = 1$, $\mathrm{CFL} = 0.1$, $p=1$ \acp{fe}. 
}
	\label{fig: advection mass conservation}
\end{figure}

Mass conservation is assessed using the advection of two cosine bells by a time-dependent, divergence free velocity field, 
as defined in Section 2.2.2 and Section 2.3 of \citet{Lauritzen2014}.  
Since the velocity field is divergence free, we use the intrinsic \ac{supg} method \cref{eq: advection supg} with $p=1$ continuous \acp{fe}, and a Crank Nicolson integrator. 
\Cref{fig: advection mass conservation} shows the initial condition is recovered at $t=5$, and that the mass conservation error reduces as the spatial refinement level increases, as expected.

\subsection{Thermal shallow water equations}
\label{sec: tsw}

A more intricate application is the \twoD thermal shallow water equations, 
which describe the evolution of fluid velocity, $\surfvec{u}$,  and depth, $\surf{\varphi}$, as a result of buoyancy transport
\cite{Ricardo2023,Ricardo2024-dg,Tambyah2025,Eldred2019,LeeRicardoTambyah2026}. 
This atmospheric system links other \twoD models that do not involve thermodynamic scalars, like the shallow water equations \cite{Cotter2012,Mcrae2013}, 
to the \threeD
compressible Euler equations typically used in operational weather models, like LFRic \cite{Adams2019}. 
The vector invariant form of the thermal shallow water equations that includes density weighted buoyancy, $\surf{B}$, transport  is
\begin{subequations}
	\label{eq: tsw}
\begin{align}
	\partial_t \surfvec{u} + 
	\surf{q} \surfvec{F}^{\perp}
	+ \surfgrad \surf{\Phi}
	+ \surf{b}\surfgrad \surf{\vartheta} = 0  
	,
	& \qquad \text{in }\thesurface \times (0,t],
	\label{eq: tsw u}
	\\
	\partial_t \surf{\varphi} + \surfdiv \surfvec{F} = 0 
	,
	& \qquad \text{in }\thesurface \times (0,t],
		\label{eq: tsw h} 
	\\
	\partial_t \surf{B} + \surfdiv (\surf{b}\surfvec{F}) = 0 
	,
	&\qquad \text{in }\thesurface \times (0,t],
		\label{eq: tsw B}
	\\
	\text{where} \quad 
	\surfvec{F} = \surf{\varphi}\surfvec{u},\quad 
	\surf{\Phi} = \frac{1}{2}\surfvec{u}\cdot \surfvec{u} + \frac{1}{2}\surf{B},\quad 
	\surf{\vartheta} = \frac{1}{2}\surf{\varphi},  \quad 
	&
	\surf{q}  = \frac{\surfgrad^{\perp} \cdot \, \surfvec{u} + \surf{\coriolis} }{\surf{\varphi}} , \quad
	\surf{b} = \frac{\surf{B}}{\surf{\varphi}},
	\label{eq: tsw constraints}
\end{align} 
\end{subequations}%
are the mass flux,
Bernoulli potential, temperature, potential vorticity, and  buoyancy, respectively, and 
$\surf{\coriolis}$ is the Coriolis parameter. 
In a previous study, the authors developed a compatible \ac{fe} \discretization 
that includes upwinded fluxes on a doubly periodic spatial domain that could represent a Cartesian domain or a curved manifold \cite{Tambyah2025}. 
We now reformulate the previous findings of \citet{Tambyah2025} 
in an intrinsic framework using the systematic derivations  in \citet{Tambyah2026_pdesphere}.

\subsubsection{Semi-discrete system}

The intrinsic formulation of the prognostics variables is: 
find $\flux{u}_h \in \V^1$, $\chart{\varphi}_h \in \V^2$, $B_h \in \V^2$, such that
\begin{subequations}
	\label{eq: tsw prognostic weak}
	\begin{align}
		\int_{\thechart} \partial_t \flux{u}_h 
		&\cdot (g \psi_{\flux{u}_h} ) \frac{1}{\sqrt{g}}\VOL 
		+ \int_{\thechart} 
		 \chart{q}_h 
		R \, \flux{F}_h
		\cdot \psi_{\flux{u}_h}  \VOL  
		- \int_{\thechart} \chart{\Phi}_h \Div \psi_{\flux{u}_h}  \VOL 
		\nonumber 
		\\
		&-a(\psi_{\flux{u}_h} ,b_h,\vartheta_h) 
		-s(\psi_{\flux{u}_h},b_h,\vartheta_h )
		= 0 
		,
		&
		\forall \psi_{\flux{u}_h} \in \V^1,
		\\ 
		\int_{\thechart}  \partial_t\chart{\varphi}_h &\phi_{\chart{\varphi}_h}\sqrt{g}\VOL
		+\int_{\thechart}\phi_{\chart{\varphi}_h} \Div  \flux{F}_h~ \VOL 
		= 0,
		&
		\forall \phi_{\chart{\varphi}_h} \in \V^2, 
		\label{eq: tsw h weak}
		\\
		\int_{\thechart} \partial_t B_h
		 \phi_{B_h} \sqrt{g}\VOL 
		&+a(\flux{F}_h,b_h,\phi_{B_h}) 
		+s(\flux{F}_h,b_h,\phi_{B_h})
		= 0  ,
	 	&
		\forall \phi_{\chart{B}_h} \in \V^2.
	\end{align} 
\end{subequations}%
The forms 
$a(\cdot,b_h,\cdot) 
$
and $
s(\cdot,b_h,\cdot) = s_{\mathrm{c}}(\cdot,b_h,\cdot)
+ s_{\mathrm{up}}(\cdot,b_h,\cdot)
$
are the \discretization of advection and stabilisation terms
\begin{subequations}
	\label{eq: tsw forms}
\begin{align}
	a(\flux{w}_h,b_h,\phi_h)
	=
	&- \int_{\thechart}  \frac{1}{2} {b}_h\flux{w}_h\cdot \grad \phi_{h}  \VOL 
	+ \int_{\thechart} \frac{1}{2} \phi_{h} {b_h}\Div  \flux{w}_h \VOL   
	\nonumber\\
	&+ \int_{\thechart} \frac{1}{2} \phi_{h} \flux{w}_h\cdot \grad b_h\VOL , 
	\\ 
	s_{\mathrm{c}}(\flux{w}_h,b_h,\phi_h)
	=
	&\int_{\chart{\mathcal{E}}_0} \frac{1}{4}
	\langle \chart{b}_h\flux{w}_h\cdot\chartcovec{n}
	\rangle
	\langle \phi_{h} \rangle \AREA
	-\int_{\chart{\mathcal{E}}_0} \frac{1}{4}
	\langle \phi_{h} \flux{w}_h\cdot\chartcovec{n}
	\rangle
	\langle \chart{b}_h\rangle \AREA,
	\\
	s_{\mathrm{up}}(\flux{w}_h,b_h,\phi_h)
	=
	&\int_{\chart{\mathcal{E}}_0} \frac{1}{2} 
	\chi(\flux{w}_h)
	\langle \phi_{h} \rangle
	\langle \chart{b}_h \rangle 
	\AREA,
\end{align} 
\end{subequations}%
for all $\flux{w}_h \in \V^1$, $\phi_h\in \V^2$, 
where $\chi(\flux{w}_h) = |\flux{w}_h\cdot\chartcovec{n}^+|/2$ is the upwinding parameter.
The approximation in \cref{eq: tsw prognostic weak,eq: tsw forms} arises from  \cref{eq: pullback,eq: integration,eq: diff op,eq: skew operators}, derived  in \citet{Tambyah2026_pdesphere}, where the intrinsic formulation of numerical fluxes follows \cref{sec: upwinding} above.
Similarly, 
the semi-discrete approximation of the diagnostics variables is:
find $\flux{F}_h\in \V^1$, $\chart{\Phi}_h \in \V^2$,
$\vartheta_h \in \V^2$,  $b_h \in \V^2$, $\chart{q}_h\in \V^0$, such that
\begin{subequations}
	\label{eq: tsw constraints weak}
	\begin{align}
		&\int_{\thechart} \flux{F}_h \cdot 
		(g \psi_{\flux{F}_h})  \frac{1}{\sqrt{g}} \VOL 
		= \int_{\thechart} \chart{\varphi}_h \flux{u}_h \cdot 
		(g \psi_{\flux{F}_h}) \frac{1}{\sqrt{g}} \VOL ,
		& \psi_{\flux{F}_h} \in \V^1, 
		\\
		&\int_{\thechart} \chart{\Phi}_h \phi_{\chart{\Phi}_h} \sqrt{g}\VOL 
		=  \int_{\thechart} \frac{1}{2} \phi_{\chart{\Phi}_h}  \flux{u}_h\cdot (g \flux{u}_h )  \frac{1}{\sqrt{g}}\VOL 
		+ \int_{\thechart} B_h   \phi_{\chart{\Phi}_h} \sqrt{g}\VOL ,
		&  \forall \phi_{\chart{\Phi}_h} \in \V^2,
		\\
		&\int_{\thechart} {\vartheta_h}\phi_{\vartheta_h} \sqrt{g}\VOL 
		= \int_{\thechart} \frac{1}{2} {\varphi_h}\phi_{\vartheta_h} \sqrt{g}\VOL  ,
		& \forall \phi_{\vartheta_h} \in \V^2,
		\\
		&\int_{\thechart} b_h\varphi_h \phi_{b_h} \sqrt{g}\VOL 
		= \int_{\thechart} B_h \phi_{b_h} \sqrt{g}\VOL ,
		& \forall \phi_{b_h}\in \V^2,
		\\
		&\int_{\thechart} \chart{q}_h \chart{\varphi}_h \xi_{\chart{q}_h} \sqrt{g}\VOL
		=   \int_{\thechart} ( R \, \flux{u}_h ) \cdot g^{-1} \grad \xi_{\chart{q}_h}  \sqrt{g} \VOL
		+ \int_{\thechart} \coriolis \xi_{\chart{q}_h} \sqrt{g} \VOL ,
		& \forall \xi_{\chart{q}_h} \in \V^0 ,  
		\label{eq: qh weak}
	\end{align}  
\end{subequations}%
where \cref{eq: qh weak} arises from using \cref{eq: skew operators} to intrinsically represent the skew divergence. 
This same weak representation of the vorticity is obtained using the skew integration by parts identity in \citet[Prop 2.10]{Tambyah2026_pdesphere}.

Consistency of the intrinsic formulation, and conservation of invariants follows from \citet{Tambyah2025}.
In the intrinsic driver, available in the 
\GridapGeoscience tutorials \cite{GridapGeoscienceTutorials}, 
the semi-discrete system is integrated using the explicit
\ac{ssprk3}
\cite{Shu1988}, where the mass systems are solved using the conjugate gradient method \cite{Manyer2024}.

\subsubsection{Results}

We now consider the modified \Galewsky test case for the thermal shallow water equations \cite{Ricardo2023,Ricardo2024-dg} to showcase the intrinsic capability of \GridapGeoscience to recover turbulent atmospheric flows that develop over long time scale.   
The original \Galewsky test case for the shallow water equations includes a mid-latitude jet that is barotropically unstable due to a Gaussian perturbation in the depth field \cite{Galewsky2004}.
The evolution of such an instability is highly sensitive to spurious waves in the vorticity field and the well-known grid imprinting error \cite{Staniforth2012}, 
which leads to a wave number 4 pattern as opposed to the accepted wave number 6 solution at day 6 \cite{Wimmer2020,Thuburn2015}. 
\citet{Kent2023} show lowest order solutions in the LFRic framework are affect by grid imprinting, 
while other studies \cite{Wimmer2020,Shipton2018} recover the accepted solution at day 6 for higher than lowest order \ac{fe} schemes.
A previous study added a perturbed buoyancy field to the original \Galewsky test case to analyse highly non-linear flows in the context of the thermal shallow water equations \cite{Ricardo2024-dg,Ricardo2023}. 
Ultimately, the \Galewsky test case \cite{Galewsky2004} yields highly non-linear flows that do not arise in the canonical Williamson test suite \cite{Williamson1992}.

\begin{figure}[h!]
	\centering 
	\includegraphics[width=0.9\textwidth]{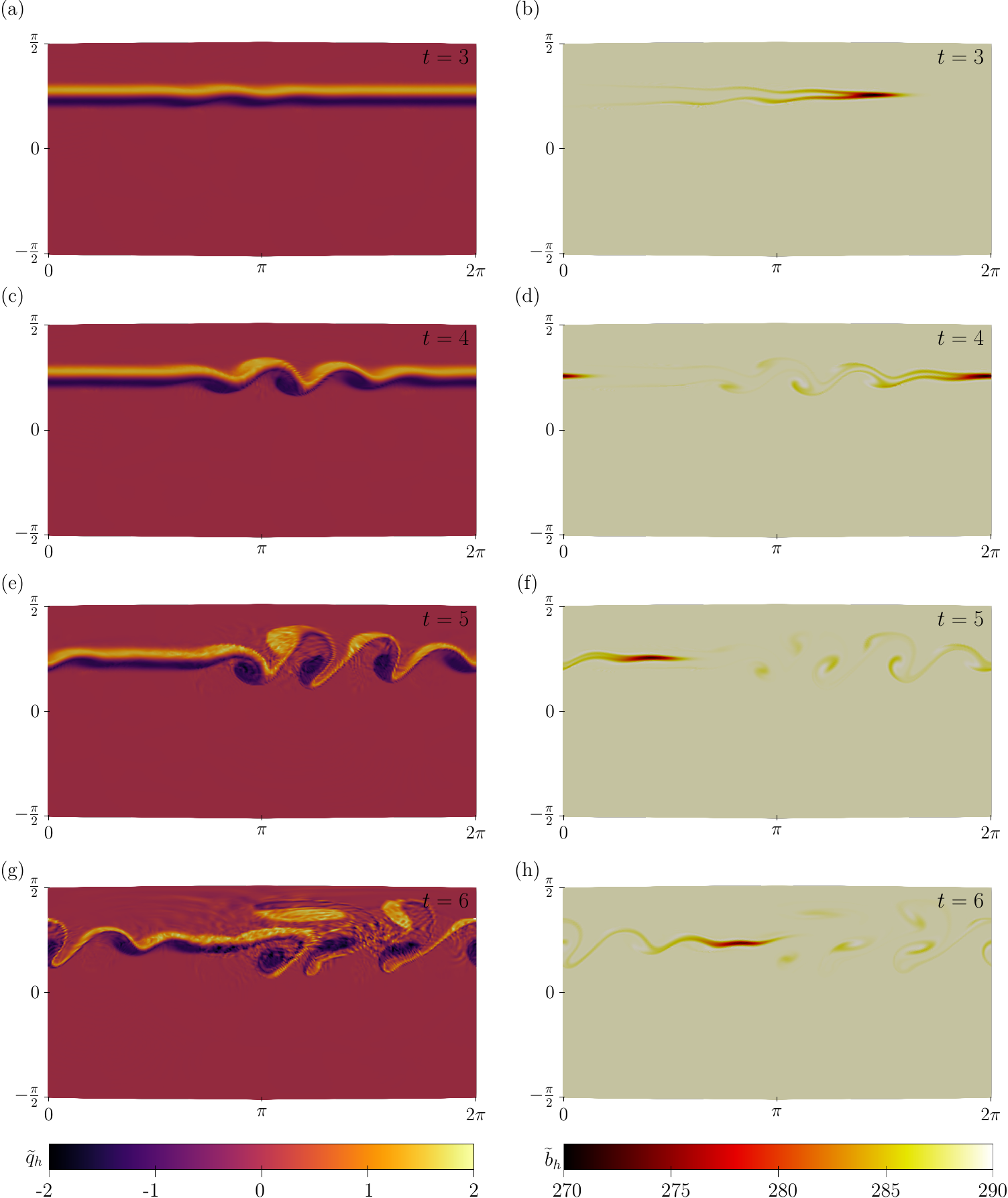}
	\caption{Modified \Galewsky test case: (a,c,e,g) show potential vorticity, and (b,d,f,h) show buoyancy, at $t = 3, 4, 5, 6$ days. 
	Parameters: $\ell=6$, $\radius = 1$,  $p=1$ \acp{fe},  CFL = 0.1.  }
	\label{fig: tsw galewsky}
\end{figure}

For the modified \Galewsky test case \cite{Ricardo2024-dg,Ricardo2023},
the initial condition is
\begin{subequations} 
	\label{eq: tsw ic}
\begin{align}
	\surfvec{u}(0)&=\surf{u}(\theta)\,\surfvec{e}_\lambda,
	\qquad
	\surf{\varphi}(0)
	= H_0 + \surf{\rho},
	\qquad
	\surf{B}(0) = \surf{b}(H_0 + \surf{\rho}),
	\\
	\text{where}  \qquad
	\surf{u}(\theta) &=
	\begin{dcases}
		\frac{u_0}{\en}  \exp\left[
	 	\{ (\theta-\theta_0 )( \theta-\theta_1 )  \}^{-1}
		\right] ,
		& \text{if } \theta_0 < \theta  < \theta_1,
		\\
		0, & \text{otherwise},
	\end{dcases}
	\\
	\surf{\rho}(\lambda,\theta) &=
	\eta \cos(\theta)\exp[  -( \mu \lambda)^2
	  -( \nu \{\theta_2 - \theta\})^2  ]
	-  \surf{\epsilon}(\theta),
	\label{eq: tsw ic rho}
	\\
	\surf{\epsilon}(\theta)
	&=\frac{\re  }{\gravity} \int_{-\pi/2}^\theta \surf{u}(\theta')
	\left[ \surf{\coriolis}(\theta') + \frac{ \surf{u}(\theta') \tan(\theta') }{\re }
	\right] \mathrm{d}\theta' ,
	\label{eq: tsw ic esp}
	\\
	\surf{b}(\lambda,\theta) &= \gravity
	\left(
	1 - 0.1 \cos(\theta)\exp\left\{ -(3\lambda)^2 -(15 [\pi/4 - \theta ]  )^2
	\right\}
	\right) ,
	\label{eq: tsw ic b}
\end{align}
\end{subequations}%
$\surf{\coriolis}(\theta) = 2 \Omega \sin (\theta)$,
and \cref{eq: tsw ic esp} is evaluated using appropriate numerical quadrature. 
This initial configuration is non-dimensionalised using the length scale $\re$ and time scale $1/\Omega$ in order to solve \cref{eq: tsw constraints weak,eq: tsw prognostic weak,eq: tsw forms} on a \twoD cubed sphere manifold with radius $\radius=1$, with  $\ell=6$ levels of refinement and  $p=1$ \acp{fe}. 
At this resolution, a quadrature degree of at least 7 is sufficient to capture the geometry of the cubed sphere to machine precision \cite[Fig. 2]{Tambyah2026_pdesphere}.  
The solution in \cref{fig: tsw galewsky} is computed with a conservative value of $q = 11$. 
The time step is $\Delta t = C h$ where $C= \mathrm{CFL}/p\sqrt{\gravity H_0}$.
The parameters are: 
$\re = 6.36 \times 10^{6}$ m, 
$\Omega = 7.292\times 10^{-5}$ s$^{-1}$,
$\gravity = 9.8$ m s$^{-2}$,
$H_0 = 10 \times 10^3$ m, 
$u_0 = 80 $ m s$^{-1}$,
$\eta = 120$ m,  
$\theta_0 = \pi/7$, $\theta_1 =\pi/2-\theta_0$, $\theta_2 = \pi/4$,
$\en = \exp[ -4/(\theta_1 - \theta_0  )^{2}]$,
$\mu = 3$, $\nu = 15$.

\cref{fig: tsw galewsky}(a,c,e,g) shows the jet evolves with the correct wave number, with the day 6 solution comparable to an accepted solution \cite[Fig. 8]{Shipton2018}.
The development of the buoyancy perturbation along the jet is shown in \cref{fig: tsw galewsky}(b,d,f,h). 
The correct wave number is recovered on a mesh coarser than those used to demonstrate the LFRic framework, which employs a $C^0$ piecewise-polynomial representation of the cubed sphere manifold \cite{Kent2023}, and other \ac{fe} studies that employ higher order methods on refined meshes \cite{Shipton2018}. We attribute this to the exact and smooth geometry of the intrinsic framework, which is free of geometric consistency errors.

\subsection{Linearised \Boussinesq equations }

As a use case of the \threeD cubed sphere manifold, 
we consider a linearisation of the compressible \Boussinesq equations. 
The continuous system is
\begin{subequations}
	\label{eq: linearised boussinesq}
	\begin{align}
		\partial_t \surfvec{u} 
		+ \omega \surfvec{k} \times \surfvec{u}
		+ \surfgrad \surf{\varphi} - \surf{b}  \surfvec{k} 
		&= 0 ,
		\quad  \text{in } \thesurface\times (0,t],
		\\
		\partial_t \surf{\varphi}  + c^2 \surfdiv \surfvec{u}
		&= 0 ,
		\quad  \text{in } \thesurface\times (0,t],
		\\
		\partial_t \surf{b} + N^2 \surfvec{u}\cdot \surfvec{k}
		&= 0,
		\quad  \text{in } \thesurface\times (0,t],
	\end{align} 
\end{subequations}%
where $\surfvec{u}$ is the fluid velocity, $\surf{\varphi}$ is the pressure, $\surf{b}$ is the buoyancy,   $\omega$ is the rotational frequency, $c$ is the speed of sound,  $N$ is the buoyancy frequency, 
$\surfvec{k}$ is the outward radial unit vector,
and  $\surfvec{u}\cdot \surfvec{k}=0$ on the radial boundary. 
In this system, $\thesurface$ represents a \threeD atmospheric shell with radius $\radius$ and thickness $\thickness$,
where  \cref{eq: 3D map} is the geometrical map. 
\cref{lst: 3D model} shows the initialisation of such a manifold in \GridapGeoscience
via the forest-of-octree data structure in p4est (lines \ref{lb3d-lstP4estStart}--\ref{lb3d-lstP4estEnd}),
which supports uniform refinement in the horizontal and vertical directions.
As an alternative, in lines \ref{lb3d-lstP6estStart}--\ref{lb3d-lstP6estEnd} we use the p6est extension of p4est  with distinct levels of horizontal and vertical refinement to obtain anisotropic refinement.
\begin{listing}[h!]
	\centering 
	\includegraphics[width=0.95\textwidth]{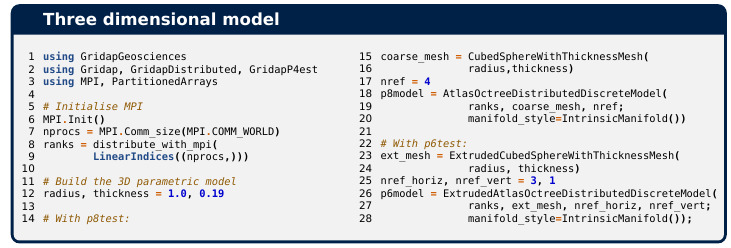}
	\caption{Initialisation of the \threeD atlas model. 
	}
	\label{lst: 3D model}
\end{listing}
\subsubsection{Semi-discrete system}

To formulate the intrinsic semi-discrete system,
we first express the radial direction parametrically. For the \threeD atmospheric shell, $n = m$ and the differential $\dmap$ is invertible. 
Consequently, the outward radial unit vector is $\surfvec{k} = \dmap \, g^{-1} \chartcovec{k} / || \dmap \, g^{-1} \chartcovec{k} ||$, which is analogous to the facet conormal \cref{eq: push normal}.
For the specific \threeD mesh illustrated in \cref{fig: cubed sphere}(c),
the local coordinates $(\alpha, \beta, \gamma)$ are defined such that $\chartcovec{k} = (0,0,1)$. 
That is, $\chartcovec{k}$ collects the covariant components of the differential of the radial coordinate $\gamma$, so that $\surfvec{k}$ points in the outward radial direction and is orthogonal to each spherical shell within $\thesurface$, which arise at constant values of $\gamma$.
The normalisation involves only the metric, $|| \dmap \, g^{-1} \chartcovec{k} ||^2 = \chartcovec{k} \cdot g^{-1} \chartcovec{k}$, so no ambient quantity enters the formulation below.
Then using \cref{eq: integration,eq: cross product,eq: diff op}
yields the semi-discrete system:
find $\flux{u}_h \in \V^1$, $\chart{\varphi}_h\in \V^2$, $\chart{b}_h\in \V^2$ such that 
\begin{subequations}
	 \label{eq: linearised boussinesq weak}
	\begin{align}  
		\int_{\thechart} \partial_t \flux{u}_h
		&\cdot (g\psi_{\flux{u}_h} ) \frac{1}{\sqrt{g}} \VOL 
		+ \int_{\thechart} \omega    
		\left\{ 
		\frac{\chartcovec{k}}{ ( \chartcovec{k} \cdot g^{-1} \chartcovec{k} )^{1/2} }
		\times 
		\left( 
		\frac{g\flux{u}_h}{\sqrt{g}}
		\right)
		\right\} \cdot  g\psi_{\flux{u}_h} 
		\frac{1}{\sqrt{g}}\VOL 
		\nonumber \\
		&- \int_{\thechart} \chart{\varphi}_h \Div \psi_{\flux{u}_h}   \VOL 
		- \int_{\thechart} \chart{b}_h  \psi_{\flux{u}_h}\cdot \frac{\chartcovec{k}}{ ( \chartcovec{k} \cdot g^{-1} \chartcovec{k} )^{1/2} } \VOL
		= 0,
		&  \forall \psi_{\flux{u}_h} \in \V^1,
		\label{eq: linearised boussinesq weak u}
		\\
		\int_{\thechart}  \partial_t\chart{\varphi}_h &\phi_{\chart{\varphi}_h} \sqrt{g}\VOL 
		+ \int_{\thechart}  c^2 \phi_{\chart{\varphi}_h} \Div  \flux{u}_h \VOL 
		= 0,
		& \forall \phi_{\chart{\varphi}_h} \in \V^2,
		\\  
		\int_{\thechart} \partial_t \chart{b}_h
		& \phi_{\chart{b}_h}  \sqrt{g}\VOL 
		+ \int_{\thechart} N^2 \phi_{\chart{b}_h}  \flux{u}_h\cdot \frac{\chartcovec{k}}{ ( \chartcovec{k} \cdot g^{-1} \chartcovec{k} )^{1/2}} \VOL 
		= 0,
	& \forall \phi_{\chart{b}_h} \in \V^2,
	\end{align} 
\end{subequations}%
where $\V^1$ is the \RT space that enforces the zero-flux boundary condition on the radial boundaries of the parametric space at $\gamma=0$ and $\gamma=1$.
Such a \ac{fe} space is initialised in the corresponding driver, available in the \GridapGeoscience tutorials \cite{GridapGeoscienceTutorials}, 
by passing the appropriate boundary tags to the \ac{fe} space constructor.

\subsubsection{Results}


Similar to other studies that also consider a linearisation of the \Boussinesq equations 	\cite{Gibson2019},
a gravity wave test is used to assess the ability of the semi-discrete formulation \cref{eq: linearised boussinesq weak} to
capture surface and vertical dynamics. 
The initial condition is
\begin{align}
	\surfvec{u}(0) &= u_0 (-y, x, 0),
	&
	\surf{\varphi}(0) &= 0,
	&
	\surf{b}(0) &= \frac{d^2}{d^2 + r_c^2} \sin \left( \frac{2\pi \eta}{ \zeta   } \right) , 
	\label{eq: LB IC}
\end{align}
where $\eta= r - \radius$ is the radial distance from the surface of the inner shell,
and $r_c= \arccos \left(
\cos(\theta)\cos(\lambda - \lambda_c)  \right) $. 
To solve \cref{eq: linearised boussinesq weak}, we use $p=1$ \acp{fe},  $\ell=4$ mesh refinement levels, 
and a Crank Nicolson time integrator
where the duration is sufficiently long for the expected buoyancy and pressure wave to develop. 
In a non-dimensional frame of reference, we choose $\radius = 1$ and $\thickness = 0.19$ to achieve the same  grid aspect ratio
as a previous study that considered the initial condition \cref{eq: LB IC} on a reduced Earth \cite{Gibson2019}. 
The non-dimensional parameters are: $u_0 = 0.058$, $\zeta = 0.38$, $d = 0.095$,
$c = 1$, $N = 1.48$, $\Omega_r = 0.01$,
$ \lambda_c= 2\pi/3$,
$\omega = 2\Omega_r \sin(\theta)$.

\begin{figure}[h!]
	\centering
	\includegraphics[width=0.9\textwidth]{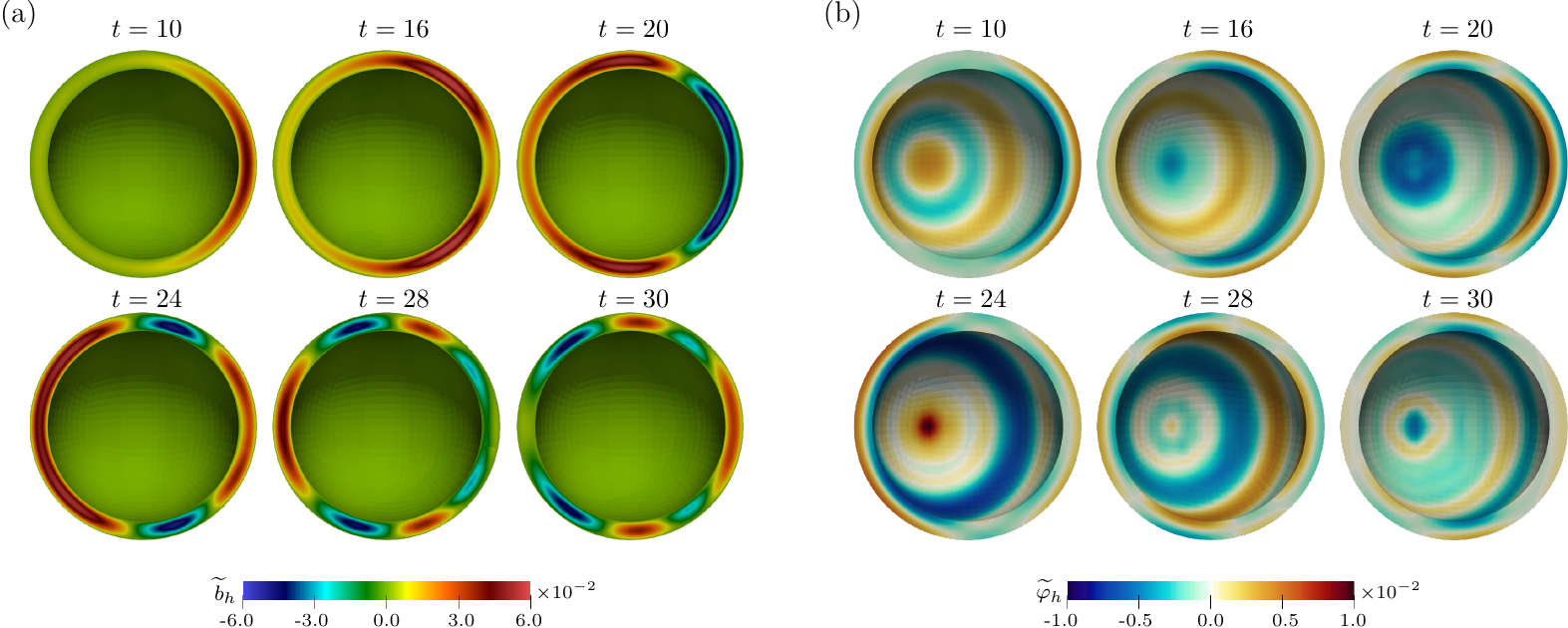}
	\caption{Evolution of (a) buoyancy and (b) pressure  for the linearised \Boussinesq equations with $p=1$ \acp{fe}. 
	Snapshots show the solution in the $x$-$z$ plane along the (0,1,0) direction, which corresponds to the view from Panel III, where time is unit-less. 
	Parameters:	$\ell=4$, $\radius=1$, $\thickness=0.19$  }
	\label{fig: LB}
\end{figure}

\Cref{fig: LB}(a) illustrates the buoyancy perturbation  on the inner- and outer-most shells is zero, 
meaning the buoyancy wave evolves in the vertical space between these boundary shells. 
Such vertical dynamics are expected since the initial buoyancy field  \cref{eq: LB IC}  depends on the radial distance from the inner sphere, 
and are only visible at a cross section of the \threeD manifold.
At the same cross section, \cref{fig: LB}(b) shows the pressure wave travels over the surface of the manifold, in contrast to the vertical buoyancy wave. 
Both the vertical buoyancy wave and surface pressure wave oscillate with time, and demonstrate the extensibility of \GridapGeoscience
to \threeD  applications.

%% file: conclusion.tex
\section{Conclusion and future work}
\label{sec: conclusion}

In this study, we develop a new method for generating discrete manifolds that relies on coarse atlas information and mesh refinement in general.  
Our abstract framework is available in \GridapGeoscience, a new Julia package that extends the \gridap \ac{fe} library to general manifolds. 
Since the high-level \gridap interface enables the seamless definition of intrinsic variational forms, 
which include general metric transformations, 
\GridapGeoscience supports both  intrinsic and extrinsic \ac{fe} formulations.
This is a key novelty of the present work in comparison to other studies that exclusively consider extrinsic functionality \cite{Rognes2013,Guba2014,FiredrakeUserManual}. 
Although not discussed in this work, 
\GridapGeoscience is a high-performance library that leverages the Julia 
\ac{jit} 
compiler, multiple dispatch and lazy evaluation of \gridap \cite{Verdugo2022,Badia2020} to provide computational efficiency. 
High resolution simulations on the cubed sphere manifold in two and three dimensions demonstrate the capacity of our implementation in \GridapGeoscience to provide intrinsic and extrinsic frameworks that successfully resolve complex flows over long periods of time.

There are many avenues to extend the functionality in \GridapGeoscience.
While the cubed sphere manifold is highlighted in this study, 
our approach to developing atlas triangulations is sufficiently general to support other manifolds. 
In particular, non-zero genus manifolds
or  manifolds that represent domains with orography via the use of terrain following coordinates \cite{Wood2013,Melvin2024_orography}. 
Since our methodology relies on mesh refinement in general, 
all meshes are nested and the development of multilevel methods on manifolds is supported \cite{Mardal2011,Benzi2005}. 
In this study, we take the most fundamental approach to mesh refinement and always consider conforming meshes that consist of elements at the same refinement level. 
The tree-based atlases developed in this study are amenable to local adaption, where mesh elements may be at different levels of refinement. 
Such non-conforming meshes are particularly advantageous in applications where the solution exhibits highly localised features \cite{Ferguson2016}. 
Extending \GridapGeoscience to non-conforming meshes introduces additional complexity in the implementation of \ac{fe} formulations \cite{Boffi2006,Arnold1985,Crouzeix1973}, 
and potentially the use of linear constraints \cite{Shephard1984}. 
Another choice we make is to demonstrate the functionality of \GridapGeoscience using space-only geometrical maps. 
One avenue of extension is to consider space-time geometrical maps and abstract \GridapGeoscience to manifolds that represent moving domains or interfaces \cite{Dziuk2013}.
This would require careful consideration of the fully discrete system and possibly space-time \ac{fe} methods \cite{BadiaDilip2023}.
Thus, the extension of \GridapGeoscience to evolving manifolds is left for future consideration.